\documentclass[12pt]{article}
\usepackage{geometry} 
\usepackage[ latin1 ]{inputenc}
\usepackage[english,french]{babel}
\usepackage{amsmath}
\usepackage{amsfonts}
\usepackage{amssymb}
\usepackage{amsthm}
\usepackage{textcomp}
\usepackage{eucal}
\usepackage{array}
\theoremstyle{plain}
\newtheorem*{theorem}{Théorème}

\theoremstyle{definition}

\theoremstyle{remark}

\title{Sur le centralisateur d'une involution de $^{2}E_{6}(2)$.}
\author{Marguerite-Marie Virotte-Ducharme\\
Équipe des groupes finis, UMR 7586;
Institut de Mathématiques\\ 175, rue du Chevaleret;
75013 Paris France.}

\begin{document}

\maketitle
\begin{abstract}
On établit que le groupe $2^{20+1}.U_{6}(2)$, centralisateur dans $^{2}E_{6}(2)$ d'une involution de la classe $2A$, est un quotient du groupe de Coxeter défini via le diagramme $Q_{222}$. Cela répond à un problème laissé ouvert dans la d\'etermination des $Q$-groupes.
\end{abstract}
\begin{otherlanguage}{english}
\begin{abstract}
In this paper we prove that $2^{20+1}.U_{6}(2)$, known as the centralizer of an involution in the group $^{2}E_{6}(2)$ is a quotient of a Coxeter group. We obtain a presentation of $2^{20+1}.U_{6}(2)$ as a $Q_{222}$-group, which now resolve a long pending question.
\end{abstract}
\end{otherlanguage}
\section{Introduction}
\subsection{Motivation et résultats}\label{01}
Rappelons qu'un groupe de Coxeter est engendré par un ensemble d'involutions correspondant aux sommets d'un graphe $X$ (graphe de Coxeter) tel que si $x$ et $y$ sont des éléments distincts dans $X$ leur produit est d'ordre $3$ s'ils sont joints par une arête et d'ordre $2$ dans le cas contraire.
On appelle $Q$-groupe un quotient du groupe de Coxeter défini par un graphe de Coxeter $Q_{rst}$ ($1\leqslant r,s,t\leqslant 4$) où $Q_{111}$ est un hexagone muni de trois bras notés $a=a_{1}, a'=a_{2},\ldots,a_{r}$, $c=c_{1},c'=c_{2},\ldots,c_{s}$, $e=e_{1},e'=e_{2},\ldots,e_{t}$:\\

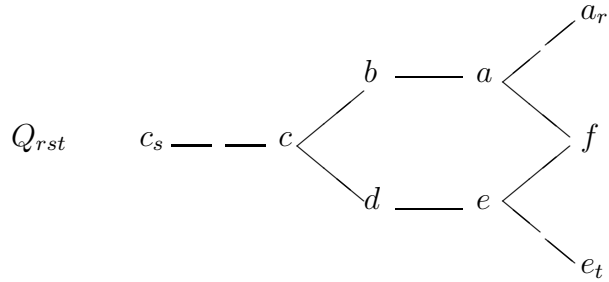
\[
\begin{picture}(150,140)

\thinlines 
\put(-32,72){$c_{s}$}
\put(-20,72){\line(1,0){15}}
\put(0,72){\line(1,0){15}}
\put(52,48){$d$}
\put(-80,72){$Q_{rst}$}
\put(20,72){$c$}
\put(64,48){\line(1,0){25}}
\put(94,48){$e$}
\put(52,96){$b$} 
\put(94,96){$a$}
\put(64,98){\line(1,0){25}}
\put(104,96){\line(6,-5){25}}
\put(104,96){\line(6,5){14}}
\put(120,110){\line(6,5){11}}
\put(133,120){$a_{r}$}
\put(104,51){\line(6,5){25}}
\put(133,72){$f$}
\put(27,72){\line(6,-5){25}}
\put(27,72){\line(6,5){25}}
\put(104,51){\line(6,-5){14}}
\put(120,37){\line(6,-5){11}}
\put(133,24){$e_{t}$}
\end{picture}
\]

            Les groupes $3^{5}\rtimes S_{6}$, $2.O_{6}^{-}(3):2$ et $U_{6}(2)$ sont des $Q$-groupes définis à partir de $Q_{111}$, $Q_{211}$ et $Q_{221}$; on en rappelle des présentations en annexe (voir aussi l'ATLAS  \cite{[2]}).\\ Dans un article déjà ancien \cite{[11]}  L. H. Soicher établit que le groupe $E=2^{3}.^{2}E_{6}(2)$ est le quotient par la relation $S=1$ du groupe de Coxeter défini par:
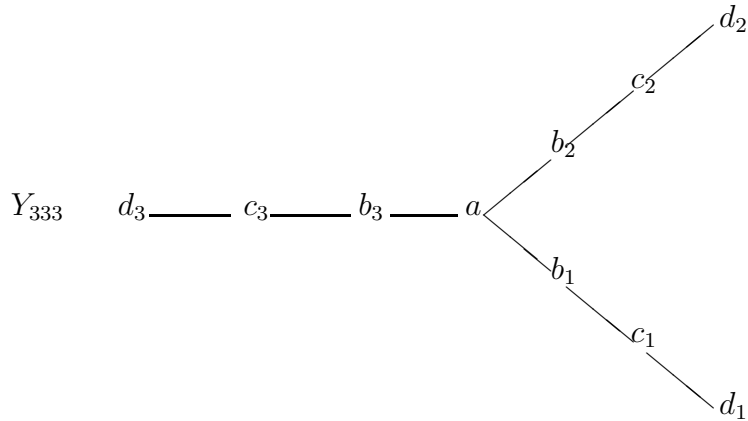
\[
\begin{picture}(150,150)
\put(-70,72){$Y_{333}$}     
\thinlines  
\put(-30,72){$d_{3}$}
\put(17,72){$c_{3}$}
\put(60,72){$b_{3}$}
\put(132,48){$b_{1}$}
\put(-18,72){\line(1,0){30}}
\put(195,-3){$d_{1}$}
\put(100,72){$a$}
\put(138,45){\line(6,-5){25}}
\put(162,24){$c_{1}$}
\put(168,20){\line(6,-5){25}}
\put(132,96){$b_{2}$} 
\put(162,120){$c_{2}$}
\put(138,98){\line(6,5){25}}
\put(195,143){$d_{2}$}
\put(168,123){\line(6,5){25}}
\put(27,72){\line(1,0){30}}
\put(72,72){\line(1,0){25}}
\put(107,72){\line(6,-5){25}}
\put(107,72){\line(6,5){25}}
\end{picture}
\]

avec $S=(ab_{1}c_{1}ab_{2}c_{2}ab_{3}c_{3})^{10}$; le groupe $E$ admet donc la présentation 
\[
E=(a,b_{i},c_{i},d_{i}(1\leqslant i \leqslant 3)/Y_{333},S=1)
\]
et son centre est $<f_{12},f_{23},f_{31}>$ où $f_{ij}=(ab_{i}b_{j}b_{k}c_{i}c_{j}d_{i})^{9}$, $\{i,j,k\}=\{1,2,3\}$ (notation de l'ATLAS \cite{[2]}).\\ L. H. Soicher met en évidence un élément $t$ de $E$ satisfaisant à
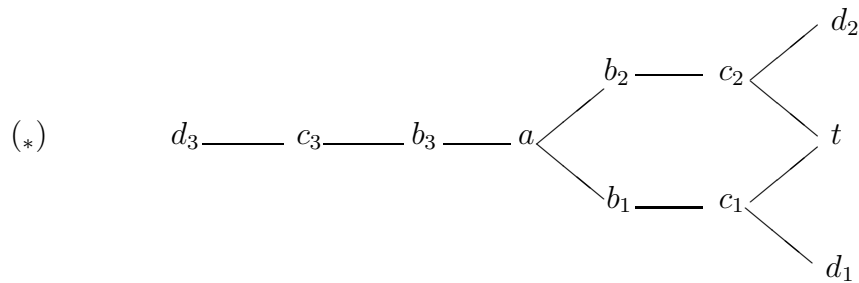
\[
\begin{picture}(150,150)
\put(-90,72){$(_*)$}    
\thinlines  
\put(-30,72){$d_{3}$}
\put(17,72){$c_{3}$}
\put(60,72){$b_{3}$}
\put(133,48){$b_{1}$}
\put(-18,72){\line(1,0){30}}
\put(215,22){$d_{1}$}
\put(100,72){$a$}
\put(144,48){\line(1,0){25}}
\put(175,48){$c_{1}$}
\put(185,48){\line(6,-5){25}}
\put(132,96){$b_{2}$} 
\put(175,96){$c_{2}$}
\put(144,98){\line(1,0){25}}
\put(217,115){$d_{2}$}
\put(187,96){\line(6,-5){25}}
\put(187,50){\line(6,5){25}}
\put(217,72){$t$}
\put(187,96){\line(6,5){25}}
\put(27,72){\line(1,0){30}}
\put(72,72){\line(1,0){25}}
\put(107,72){\line(6,-5){25}}
\put(107,72){\line(6,5){25}}
\end{picture}
\]
et à la relation hexagonale: $V=1$ avec $V=(atb_{2}c_{1}c_{2}b_{1})^{4}$.\\
Le centralisateur de $d_{3}$ dans $E$ contient donc des éléments satisfaisant aux relations de graphe $Q_{222}$ obtenu à partir de $(_*)$ en omettant $d_{3}$ et $c_{3}$.\\
On connaît la structure du centralisateur d'une involution $d$ de $E$ provenant de la classe $2A$ dans $\bar{E} =E/Z(E)$: 
$C_{E}(d)\simeq 2^{4}.2^{20}.U_{6}(2)$ et $C_{\bar{E}}(\bar{d})\simeq 2^{20+1}.U_{6}(2)$, \\$\bar{d}$ désignant l'image de $d$ dans $\bar{E}$ (\cite{[2]}, \cite{[6]}, \cite{[9]}, \cite{[12]}).\\

Dans \cite{[11]} L. H. Soicher ne donne pas de présentation de $C_{E}(d)$, mais conjecture, comme dans l'ATLAS (\cite{[2]}, \cite{[10]}) que $C_{E}(d)$ est un quotient du groupe défini par $Q_{222}$ et $V=1$.\\

L'objet de ce travail est de décrire un jeu de relations noté $rel(i)$ ($1\leqslant i \leqslant 3$) (avec $rel(1) \Rightarrow rel(2) \Rightarrow rel(3)$) suffisant pour définir les groupes $2^{20+1}.U_{6}(2)$, $2^{4}.2^{20}.U_{6}(2)$ et un troisième groupe $2^{2}.2^{4}.2^{20}.U_{6}(2)$.\\
On se propose d'établir le résultat suivant:\\
\begin{theorem}
Soit $G$ un groupe avec la présentation:
\[
G=(a,b,\ldots ,f,a',c',e'/Q_{222},V=1,rel(i))(1\leqslant i \leqslant 3) 
\]
où
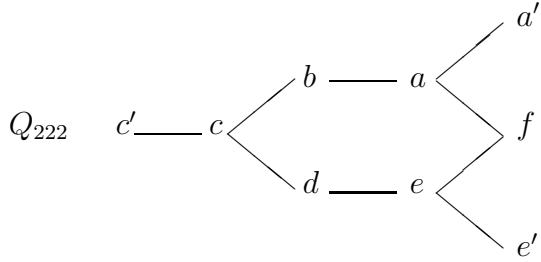
\[
\begin{picture}(150,150)
\put(-70,72){$Q_{222}$} 
\thinlines  
\put(-30,72){$c'$}
\put(-23,72){\line(1,0){25}}
\put(5,72){$c$}
\put(12,72){\line(6,-5){25}}
\put(12,72){\line(6,5){25}}
\put(40,90){$b$}
\put(40,50){$d$}
\put(50,92){\line(1,0){25}}
\put(50,50){\line(1,0){25}}
\put(80,90){$a$}
\put(80,50){$e$}
\put(90,93){\line(6,5){25}}
\put(90,50){\line(6,-5){25}}
\put(120,113){$a'$} 
\put(120,26){$e'$}
\put(90,50){\line(6,5){25}}
\put(90,92){\line(6,-5){25}}
\put(120,72){$f$}
\end{picture}
\]
\begin{eqnarray*}
V=(adbecf)^{4}&\mbox{(relation hexagonale)}\\
rel(i)&\mbox{(voir ci-après),}
\end{eqnarray*}
et soit $H$ le sous-groupe de $G$ engendré par $a,b,\ldots ,f,a',c'$.\\
Pour $i=1$, le groupe $G$ est isomorphe au centralisateur d'une involution de la classe $2A$ de $^{2}E_{6}(2)$: on a $G \simeq 2^{20+1}.U_{6}(2)$ et $H$ est isomorphe à $2.U_{6}(2)$.\\
Pour $i=2$, le groupe $G$ est isomorphe à $2^{4}.2^{20}.U_{6}(2)$ et $H$ à $2.2.U_{6}(2)$.\\
Pour $i=3$, le groupe $G$ est isomorphe à $2^{2}.2^{4}.2^{20}.U_{6}(2)$ et $H$ à $2^{2}.2.U_{6}(2)$.\\
Notations:
\begin{itemize}
  \item $rel(1)=\{r=1,z_{1}=z_{2}=z_{3}=1,m_{a}=m_{e}=1\}$,
  \item $rel(2)=\{r=1,\mathcal{R}=1,m_{a}=m_{e}=1\}$,
 \item $rel(3)=\{r=1,\mathcal{R}=1\}$,

 \item $r=(\mathcal{C}^{ee'de}\mathcal{A})^{4}$ avec $\mathcal{C}=c^{c'dbc}$ et $\mathcal{A}=a^{a'bfa}$,
\item $\mathcal{R}=(\mathcal{C}^{ee'de}.\mathcal{C}^{aa'ba})^{2}$,
 \item $z_{1}=(cc'bdaee')^{9}$,
  \item $z_{2}=(cc'bdeaa')^{9}$,
  \item $z_{3}=(ee'dfcaa')^{9}$,
\item $m_{a}=(aa'fbcd)^{5}$,
\item $m_{e}=(ee'dfab)^{5}$.
\end{itemize}

\end{theorem}
\subsection{Méthodes et Plan}
Notons $E$ un groupe isomorphe à $2^{3}.^{2}E_{6}(2)$ et $G$ un groupe satisfaisant aux conditions du théorème dans le cas $i=3$.\\
Le travail comporte quatre parties: des résultats préliminaires, la construction d'un certain sous-groupe $N$, la démonstration proprement dite du théorème enfin des tables et une annexe.
\subsubsection{Les résultats préliminaires}
D'abord on donne quelques compléments sur les éléments intervenant dans les relations $rel(i)$ $(1\leqslant i \leqslant 3)$. Puis on vérifie que les relations qui déterminent $G$ sont satisfaites dans $E$ de sorte que l'on a un morphisme $\Phi$ de $G$ dans $E$; $\Phi(G)$ est un sous-groupe de $E$ centralisé par une involution du système générateur de $E$.
\subsubsection{Le sous-groupe $N$}\label{N}
Sa construction est une partie essentielle de la démonstration. En remplaçant le générateur $e'$ par un de ses conjugués $a^{o}$ dans $G$, on obtient un nouveau système générateur de $G$ et on montre, grâce à la relation $R=1$, que $a'a^{o}$ est d'ordre 2. La fermeture normale de $a'a^{o}$ dans $G$ est alors le sous-groupe $N$ en question. Pour construire $N$, on détermine des conjugués de $a'a^{o}$ de manière à obtenir un système générateur $\Gamma$ de $N$. Les calculs sont longs et fastidieux, ils ne sont pas reproduits ici; certains sont rassemblés sous forme de tables: table des conjugués $\gamma^{y},\gamma \in \Gamma,y \in \{a,b,\ldots,f,a',c'\}$, table des commutateurs des éléments de $\Gamma$ \ldots. On détermine l'ordre de $N$, ses éléments centraux et son groupe des commutateurs.
\subsubsection{La preuve du théorème}
Dans chacune des situations $i=1,2,3$, on détermine les éléments centraux de $G$ et donc l'ordre de $Z(G)$. On établit les résultats concernant le sous-groupe $H$ (du théorème) puis on vérifie que $G/Z(G).N$ est isomorphe à $U_{6}(2)$; enfin on conclut que l'image de $G$ dans $E$ est bien le centralisateur d'une certaine involution de $E$.
\subsubsection{Tables. Annexe}
Pour faciliter la lecture, nous avons choisi de ne pas donner les détails des démonstrations. Nous donnons sous forme de tables des résultats concernant le groupe $N$ (\ref{N}). En annexe nous avons rassemblé des compléments utiles à la clarté du texte. En général ces résultats sont connus ou se démontrent sans difficultés; ils concernent les groupes $3^{5}\rtimes S_{6}$, $2.O^{-}_{6}(3)$, $W(E_{7})$ et $U_{6}(2)$, groupes qui admettent des présentations via les graphes de Coxeter $Q_{111}$, $Q_{211}$, $Y_{321}$ et $Q_{222}$.
\section{Préliminaires}
Dans cette section $E$ et $G$ désignent respectivement des groupes donnés avec leur présentation:
\[
E=(a,b_{i},c_{i},d_{i}(1\leqslant i\leqslant 3)/Y_{333},S=1)
\]
\[
G=(a,b,\ldots,f,a',c',e'/Q_{222},V=1,rel(3)).
\]
\subsection{Résultats concernant le groupe $E$}\label{E}
\subsubsection{Les éléments $f_{ij}$}
Pour $\{i,j,k\}=\{1,2,3\}$, on note $W_{ij}$ le sous-groupe de $E$ 
\[
W_{ij}=<a,b_{i},c_{i},d_{i},\\b_{j},c_{j},b_{k}>;
\]
 $W_{ij}$ est centralisé par $d_{k}$ et admet un élément supplémentaire $t_{ij}$ tel que:
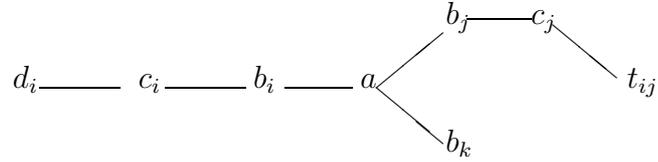
\[
\begin{picture}(150,120)
\thinlines 
\put(-30,72){$d_{i}$}
\put(17,72){$c_{i}$}
\put(60,72){$b_{i}$}
\put(132,48){$b_{k}$}
\put(-20,72){\line(1,0){30}}
\put(100,72){$a$}
\put(132,96){$b_{j}$} 
\put(164,96){$c_{j}$}
\put(140,98){\line(1,0){24}}
\put(172,96){\line(6,-5){24}}
\put(200,72){$t_{ij}$}
\put(27,72){\line(1,0){30}}
\put(72,72){\line(1,0){25}}
\put(106,72){\line(6,-5){25}}
\put(106,72){\line(6,5){25}}
\end{picture}
\]
Rappelons que $W_{ij}$ et $W_{ik}$ engendrent un sous-groupe $O_{i}$ de $E$ isomorphe à $O_{7}(3)\times 2$ dont l'involution centrale est $f_{ij}=f_{ik}$ ($f_{ij}=(ab_{i}b_{j}b_{k}c_{i}c_{j}d_{i})^{9}$) (notations de l'ATLAS) et que la relation $S=1$ impose $t_{ij}=t_{ik}$ ([2], [3], [11], [16]).\\
On désigne par $t$ l'élément $t_{31}=t_{32}$ et par $f_{i}$ l'involution $f_{ij}=f_{ik}$; $f_{i}$ est centrale dans $E$ et l'on a $Z(E)=<f_{12},f_{23},f_{31}>$ ([11]). En outre on a $f_{3}=d_{3}m$ avec $m=(ab_{1}b_{2}b_{3}c_{2}t)^{5}$ (voir Annexe 3).
\subsubsection{La relation hexagonale}\label{H}
Soit $H_{0}$ le sous-groupe de $E$ engendré par $a,b_{1},c_{1},t,b_{2},c_{2}$; ces éléments satisfont aux relations ci-dessous:
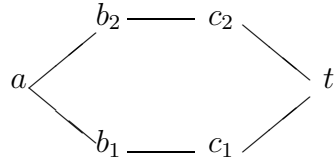
\[
\begin{picture}(150,120)
\thinlines             
\put(52,48){$b_{1}$}
 \put(20,72){$a$}
 \put(64,48){\line(1,0){25}}
 \put(94,48){$c_{1}$}
 \put(52,96){$b_{2}$} 
 \put(94,96){$c_{2}$}
 \put(64,98){\line(1,0){25}}
 \put(107,96){\line(6,-5){25}}
 \put(107,48){\line(6,5){25}}
 \put(137,72){$t$}
 \put(27,72){\line(6,-5){25}}
 \put(27,72){\line(6,5){25}}
 \end{picture}
\]
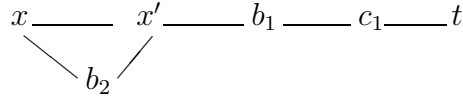
\[
\begin{picture}(150,60)
\thinlines  
\put(-30,72){$x$}
\put(17,72){$x'$}
\put(60,72){$b_{1}$}
\put(-22,72){\line(1,0){30}}
\put(100,72){$c_{1}$}
\put(110,72){\line(1,0){23}}
\put(-25,68){\line(5,-4){20}}
\put(-2,48){$b_{2}$} 
\put(10,52){\line(5,6){13}}
\put(135,72){$t$}
\put(27,72){\line(1,0){30}}
\put(72,72){\line(1,-0){25}}
\end{picture}
\]
où $x=a^{b_{1}c_{1}tc_{2}}$ et $x'=b_{1}^{c_{1}tc_{2}}$.\\
La relation $S=1$ impose que le sous-groupe $H_{0}$ soit isomorphe à $H_{3,6}$ ou à $H_{3,6}/Z(H_{3,6})$ (voir l'Annexe 1). Or $H_{0}$ est contenu dans $O_{i}$ ($O_{i} \simeq O_{7}(3)\times 2$) et son centre $Z(H_{0})$ (qui est un 3-groupe) est central dans $O_{i}$, il s'ensuit que $Z(H_{0})=1$; on a alors $(atb_{1}c_{2}c_{1}b_{2})^{4}=1$ (relation hexagonale) (Annexe 1, \cite{[8]}, \cite{VD4}, \cite{Z1}).
\subsection{Résultats concernant le groupe $G$}
\subsubsection{Notations}\label{not}
On pose $\mathcal{A}=a^{a'bfa}$,  $\mathcal{C}=c^{c'dbc}$ et $\mathcal{E}=e^{e'dfe}$; la relation $V=1$ impose que $\mathcal{A}$, $\mathcal{C}$ et $\mathcal{E}$ commutent deux à deux et que l'on a:
\[
\mathcal{A}^{edfe}=\mathcal{A}^{cbdc}, \mathcal{C}^{abfa}=\mathcal{C}^{edfe}, \mathcal{E}^{cbdc}=\mathcal{E}^{abfa}.
\]
\subsubsection{Les éléments $z_{i}$ ($1\leqslant i\leqslant 3$)}\label{zi}
Ils désignent les involutions centrales des groupes $W_{i}$  ($1\leqslant i\leqslant 3$) isomorphes à $W(E_{7})$ définis respectivement à partir de:
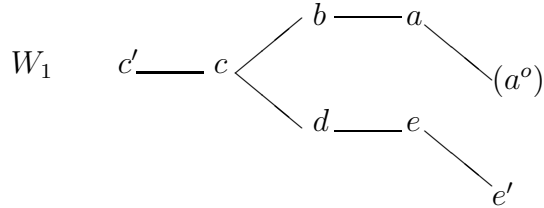
\[
\begin{picture}(150,110)
\thinlines  
\put(-30,72){$c'$}
\put(-23,72){\line(1,0){25}}
\put(6,72){$c$}
\put(14,72){\line(6,-5){25}}
\put(14,72){\line(6,5){25}}
\put(43,90){$b$}
\put(43,50){$d$}
\put(51,50){\line(1,0){25}}
\put(51,93){\line(1,0){25}}
\put(78,50){$e$}
\put(78,90){$a$}
\put(110,66){$(a^{o})$}
\put(85,90){\line(6,-5){25}}
\put(-70,72){$W_{1}$}                
\put(110,22){$e'$}
\put(85,50){\line(6,-5){25}}
\end{picture}
\]
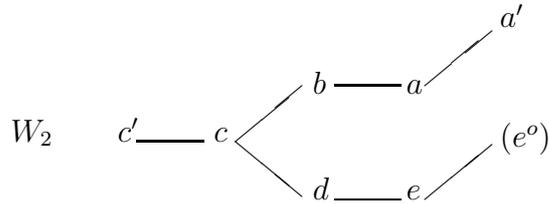
\[
\begin{picture}(150,120)
\thinlines  
\put(-30,72){$c'$}
\put(-23,72){\line(1,0){25}}
\put(6,72){$c$}
\put(14,72){\line(6,-5){25}}
\put(14,72){\line(6,5){25}}
\put(43,90){$b$}
\put(43,50){$d$}
\put(51,50){\line(1,0){25}}
\put(51,93){\line(1,0){25}}
\put(78,50){$e$}
\put(78,90){$a$}
\put(85,93){\line(6,5){25}}
\put(113,115){$a'$} 
\put(85,50){\line(6,5){25}}
\put(113,70){$(e^{o})$}
\put(-70,72){$W_{2}$}                
\end{picture}
\]
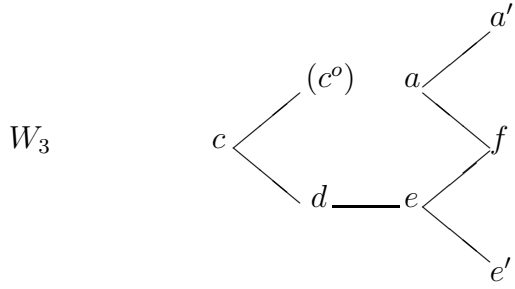
\[
\begin{picture}(150,120)
\thinlines
\put(6,72){$c$}
\put(14,72){\line(6,-5){25}}
\put(43,50){$d$}
\put(41,95){$(c^{o})$}
\put(78,50){$e$}
\put(78,95){$a$}
\put(110,72){$f$}
\put(110,72){\line(-6,5){25}}
\put(85,95){\line(6,5){25}}
\put(-70,72){$W_{3}$}                
\put(110,118){$a'$} 
\put(110,22){$e'$}
\put(85,50){\line(6,5){25}}
\put(85,50){\line(6,-5){25}}
\put(51,50){\line(1,0){25}}
\put(14,72){\line(6,5){25}}
\end{picture}
\]
On note $a^{o}$, $e^{o}$, $c^{o}$ les involutions qui permettent d'obtenir les diagrammes complétés. On a $a^{o}=\mathcal{C}^{abedcc'e'edcba}$ et $e'=\mathcal{C}^{abedcc'a^{o}abcde}$.
Les éléments centraux $z_{i}$ s'écrivent comme produit de sept involutions commutant deux à deux; (avec les notations \ref{not}) on a:
\begin{itemize}
\item $z_{1}=c'bd\mathcal{C}\mathcal{C}^{ee'de}a^{o}e'=(cc'bdaee')^{9}$
\item $z_{2}=c'bd\mathcal{C}\mathcal{C}^{aa'ba}e^{o}a'=(cc'bdeaa')^{9}$
\item $z_{3}=e'df\mathcal{E}\mathcal{E}^{aa'fa}c^{o}a'=(ee'dfcaa')^{9}$.
\end{itemize}
(voir aussi annexe 3).
\subsubsection{La relation $R=1$}\label{123}
Les relations $R=1$ et $(a'a^{o})^{2}=1$ sont équivalentes. On a $R=1$ si et seulement si l'un des $z_{i}$ est dans $Z(G)$. Sous l'hypothèse $R=1$, les éléments $z_{1}$,  $z_{2}$ et  $z_{3}$ sont centraux dans $G$.
\subsubsection{Les éléments $m_{a}$ et $m_{e}$}\label{Ra}
On désigne par $R_{a}$, $R_{c}$ et $R_{e}$ les sous-groupes de $G$ respectivement  définis à partir de:
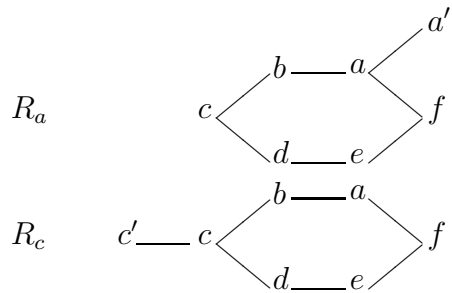
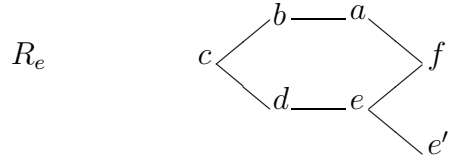
\[
\begin{picture}(100,140)
\thinlines  
\put(0,72){$c$}
\put(28,55){$d$}
\put(28,87){$b$}
\put(64,89){\line(6,-5){20}}
\put(57,55){$e$}
\put(57,89){$a$}
\put(86,72){$f$}
\put(64,89){\line(6,5){20}}
\put(-70,72){$R_{a}$}                
\put(86,105){$a'$} 
\put(64,55){\line(6,5){20}}
\put(35,89){\line(1,0){20}}
\put(35,55){\line(1,0){20}}
\put(7,72){\line(6,-5){20}}
\put(7,72){\line(6,5){20}}
\end{picture}
\] 
\[
\begin{picture}(100,40)
\thinlines  
\put(-30,72){$c'$}
\put(0,72){$c$}
\put(28,55){$d$}
\put(28,87){$b$}
\put(-3,72){\line(-20,0){20}}
\put(57,55){$e$}
\put(57,89){$a$}
\put(86,72){$f$}
\put(64,89){\line(6,-5){20}}
\put(-70,72){$R_{c}$}                
\put(64,55){\line(6,5){20}}
\put(35,89){\line(1,0){20}}
\put(35,55){\line(1,0){20}}
\put(7,72){\line(6,-5){20}}
\put(7,72){\line(6,5){20}}
\end{picture}
\]
\[
\begin{picture}(100,40)
\thinlines                    
\put(0,72){$c$}
\put(28,55){$d$}
\put(28,87){$b$}
\put(57,55){$e$}
\put(57,89){$a$}
\put(86,72){$f$}
\put(64,89){\line(6,-5){20}}
\put(-70,72){$R_{e}$}                
\put(86,38){$e'$}
\put(64,55){\line(6,5){20}}
\put(64,55){\line(6,-5){20}}
\put(35,89){\line(1,0){20}}
\put(35,55){\line(1,0){20}}
\put(7,72){\line(6,-5){20}}
\put(7,72){\line(6,5){20}}
\end{picture}
\]
Pour chacun d'entre eux la relation hexagonale est satisfaite.\\

Ces groupes sont isomorphes à $2.O^{-}_{6}(3):2$, leur involution centrale est un produit de six involutions commutant deux à deux, on la note respectivement:
\begin{itemize}
\item $m_{a}=a'bfd\mathcal{A}\mathcal{A}^{cbdc}=(aa'fbcd)^{5}$
\item $m_{c}=c'bdf\mathcal{C}\mathcal{C}^{edfe}=(cc'bdef)^{5}$
\item $m_{e}=e'dfb\mathcal{E}\mathcal{E}^{abfa}=(ee'dfab)^{5}$
\end{itemize}
(voir Annexe 3).
\subsubsection{Le sous-groupe $K_{e}$}\label{Ke}
Soit $K_{e}$ le sous-groupe de $G$ engendré par $a,b,c,d,e,f,a',c'$; $K_{e}$ est isomorphe à $2^{2}.2.U_{6}(2)$ et l'on a $Z(K_{e})=<m_{a},m_{c}>.<m_{c}z_{2}>$ (\ref{zi}, \ref{123}). Soit $e^{o}$ l'unique élément de $K_{e}$ tel que l'on ait:
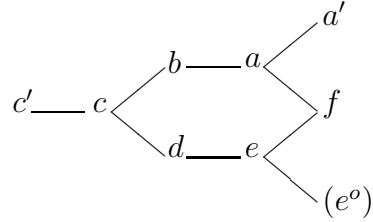
\[
\begin{picture}(100,120)
\thinlines  
\put(-30,72){$c'$}
\put(0,72){$c$}
\put(28,55){$d$}
\put(28,87){$b$}
\put(-23,72){\line(1,0){20}}
\put(57,55){$e$}
\put(57,89){$a$}
\put(86,72){$f$}
\put(64,89){\line(6,5){20}}
\put(64,89){\line(6,-5){20}}
\put(86,105){$a'$} 
\put(86,35){$(e^{o})$}
\put(64,55){\line(6,5){20}}
\put(64,55){\line(6,-5){20}}
\put(35,89){\line(1,0){20}}
\put(35,55){\line(1,0){20}}
\put(7,72){\line(6,-5){20}}
\put(7,72){\line(6,5){20}}
\end{picture}
\]
On a alors $m_{a}m_{c}=n_{e}$ où $n_{e}$ est l'involution centrale du sous-groupe\\ $<b,c,d,e,f,e^{o}>$(voir Annexe 4). On a des résultats similaires pour les sous-groupes $K_{a}=<a,b,c,d,e,f,c',e'>$ et $K_{c}=<a,b,c,d,e,f,a',e'>$\\ (voir le tableau 2.2.6).
\subsubsection{Un grand tableau}\label{126}
(Avec les notations ci-dessus). Pour chacun des groupes $K$ ($K=K_{a},K_{c},K_{e})$, le tableau ci-dessous indique un système générateur de $K$ et de trois sous-groupes $W$ isomorphes à $W(E_{7})$, $W/Z(W)$ étant un représentant de chacune des trois classes de $W^{*}(E_{7})$ de $K/Z(K)$ (\cite{VD1}, \cite{VD2}). Pour chaque groupe, on précise les éléments centraux (tableau voir page suivante).
\subsubsection{Les éléments $m_{a}$, $m_{c}$ et $m_{e}$}\label{127}
Les éléments $m_{a}$, $m_{c}$ et $m_{e}$ sont centraux dans $G$. En effet $m_{a}$ et $m_{c}$ (resp. $m_{c}$ et $m_{e}$) appartiennent à $Z(K_{e})$ (resp. $Z(K_{a})$) et ils centralisent $e'$ (resp. $a'$). (\ref{Ra}, \ref{Ke}).
\pagebreak
\[
\begin{array}{|p{55 pt}| >{$}c<{$} | >{$}c<{$} | >{$}c<{$} |}
\hline

 & $K_{a}$ &  $K_{e}$ & $K_{c}$\\ \hline
\raisebox{72 pt}{\parbox {55 pt}{ Système générateur de $K$} }
&
\begin{picture}(100,120)
\thinlines  
\put(-5,72){$c'$}
\put(20,72){$c$}
\put(39,57){$d$}
\put(39,83){$b$}
\put(4,72){\line(1,0){12}}
\put(58,57){$e$}
\put(58,83){$a$}
\put(78,72){$f$}
\put(67,80){\line(6,-5){10}}
\put(67,61){\line(6,5){10}}
\put(46,83){\line(1,0){10}}
\put(46,58){\line(1,0){10}}
\put(27,72){\line(6,-5){10}}
\put(27,72){\line(6,5){10}}
\put(67,57){\line(6,-5){10}}
\put(78,45){$e'$}
\end{picture}

 &
\begin{picture}(100,120)
\thinlines  
\put(20,72){$c$}
\put(-5,72){$c'$}
\put(4,72){\line(1,0){12}}
\put(39,57){$d$}
\put(39,83){$b$}
\put(67,80){\line(6,-5){10}}
\put(67,61){\line(6,5){10}}
\put(58,57){$e$}
\put(58,83){$a$}
\put(78,72){$f$}
\put(67,83){\line(6,5){10}}
\put(78,93){$a'$} 
\put(46,83){\line(1,0){10}}
\put(46,58){\line(1,0){10}}
\put(27,72){\line(6,-5){10}}
\put(27,72){\line(6,5){10}}
\end{picture}
 
 &
\begin{picture}(100,120)
\thinlines                     
\put(0,72){$c$}
\put(19,57){$d$}
\put(18,83){$b$}
\put(38,57){$e$}
\put(38,83){$a$}
\put(61,72){$f$}
\put(47,80){\line(6,-5){10}}
\put(47,57){\line(6,5){10}}
\put(61,45){$e'$}
\put(47,83){\line(6,5){10}}
\put(47,57){\line(6,-5){10}}
\put(26,83){\line(1,0){10}}
\put(26,58){\line(1,0){10}}
\put(7,72){\line(6,-5){10}}
\put(7,72){\line(6,5){10}}
\put(61,95){$a'$} 
\end{picture}
\\ \hline
\parbox{55 pt}{Système générateur de W, involution centrale} 
&\begin{picture}(100,120)
\thinlines  
\put(-5,72){$c'$}
\put(4,72){\line(1,0){10}}
\put(18,72){$c$}
\put(25,72){\line(6,-5){10}}
\put(25,72){\line(6,5){10}}
\put(38,80){$b$}
\put(38,58){$d$}
\put(45,58){\line(1,0){10}}
\put(45,80){\line(1,0){10}}
\put(58,58){$e$}
\put(58,80){$a$}
\put(79,72){$(a^{o})$}
\put(67,80){\line(6,-5){10}}
\put(79,42){$e'$}
\put(67,58){\line(6,-5){10}}
\put(95,58){$z_{1}$}
\end{picture}

& \begin{picture}(100,120)
\thinlines 
\put(-5,72){$c'$}
\put(4,72){\line(1,0){10}}
\put(18,72){$c$}
\put(25,72){\line(6,-5){10}}
\put(25,72){\line(6,5){10}}
\put(38,80){$b$}
\put(38,58){$d$}
\put(45,58){\line(1,0){10}}
\put(45,80){\line(1,0){10}}
\put(58,58){$e$}
\put(58,80){$a$}
\put(65,80){\line(6,5){10}}
\put(76,90){$a'$} 
\put(65,58){\line(6,-5){10}}
\put(76,50){$(e_{o})$}
\put(90,72){$z_{2}$}
\end{picture}

& \begin{picture}(100,120)
\thinlines 
          \put(0,72){$c$}
         \put(19,57){$d$}
         \put(17,83){$(c^{o})$}
         \put(38,57){$e$}
         \put(38,83){$a$}
         \put(61,72){$f$}
         \put(47,80){\line(6,-5){10}}
         \put(47,57){\line(6,5){10}}
         \put(61,45){$e'$}
         \put(47,83){\line(6,5){10}}
         \put(47,57){\line(6,-5){10}}
         \put(26,58){\line(1,0){10}}
         \put(7,72){\line(6,-5){10}}
         \put(7,72){\line(6,5){10}}
         \put(61,95){$a'$}
         \put(80,72){$z_{3}$}
         \end{picture}\\
 &\begin{picture}(100,40)
\thinlines  
        \put(-5,72){$c'$}
        \put(4,72){\line(1,0){10}}
        \put(18,72){$c$}
        \put(25,72){\line(6,-5){10}}
        \put(67,80){\line(6,-5){10}}
        \put(78,72){$f$}
        \put(38,58){$d$}
        \put(45,58){\line(1,0){10}}
        \put(67,58){\line(6,5){10}}
        \put(58,58){$e$}
        \put(58,80){$a$}
        \put(78,90){$(a^{o})$}
        \put(67,80){\line(6,5){10}}
        \put(95,72){$z'_{1}$}                
        \put(79,42){$e'$}
        \put(67,58){\line(6,-5){10}}
        \end{picture}
            & \begin{picture}(100,40)
\thinlines  
       \put(-5,72){$c'$}
       \put(4,72){\line(1,0){10}}
       \put(18,72){$c$}
       \put(67,80){\line(6,-5){10}}
       \put(25,72){\line(6,5){10}}
       \put(38,80){$b$}
       \put(78,72){$f$}
       \put(67,58){\line(6,5){10}}
       \put(45,80){\line(1,0){10}}
       \put(58,58){$e$}
       \put(58,80){$a$}
       \put(65,80){\line(6,5){10}}
       \put(76,90){$a'$} 
       \put(65,58){\line(6,-5){10}}
      \put(78,50){$(e_{o})$}
       \put(95,72){$z'_{2}$}                
       \end{picture}
&\begin{picture}(100,40)
\thinlines                    
	 \put(0,72){$c$}
          \put(12,100){$(c^{o})$}
          \put(7,72){\line(6,5){10}}
          \put(30,80){\line(1,0){10}}
          \put(20,80){$b$}
         \put(38,57){$e$}
         \put(38,83){$a$}
         \put(61,72){$f$}
         \put(47,80){\line(6,-5){10}}
         \put(47,57){\line(6,5){10}}
         \put(61,45){$e'$}
         \put(47,83){\line(6,5){10}}
         \put(47,57){\line(6,-5){10}}
         \put(7,78){\line(1,5){4}}
         \put(61,95){$a'$} 
         \put(85,72){$z'_{3}$}
         \end{picture}\\
         &\begin{picture}(100,80)
\thinlines 
         \put(-5,72){$c'$}
         \put(4,72){\line(1,0){10}}
         \put(18,72){$c$}
         \put(25,72){\line(6,5){10}}
         \put(38,80){$b$}
         \put(45,80){\line(1,0){10}}
         \put(67,80){\line(6,-5){10}}
         \put(78,72){$f$}
         \put(67,58){\line(6,5){10}}
         \put(58,58){$e$}
         \put(58,80){$a$}
         \put(78,90){$a^{o}$}
         \put(67,80){\line(6,5){10}}
        \put(78,48){$(e')$}
        \put(95,72){$z''_{1}$}                
        \put(67,58){\line(6,-5){10}}
        \end{picture}
&\begin{picture}(100,80)
\thinlines  
       \put(-8,72){$(c')$}
           \put(8,72){\line(1,0){8}}
            \put(18,72){$c$}
            \put(25,72){\line(6,-5){10}}
            \put(38,58){$d$}
            \put(45,58){\line(1,0){10}}
            \put(67,80){\line(6,-5){10}}
            \put(78,72){$f$}
           \put(67,58){\line(6,5){10}}
            \put(58,58){$e$}
             \put(58,80){$a$}
             \put(65,80){\line(6,5){10}}
             \put(76,90){$a'$} 
             \put(65,58){\line(6,-5){10}}
             \put(76,45){$e_{o}$}
              \put(90,72){$z''_{2}$}                
       	 \end{picture}
	 &\begin{picture}(100,80)
\thinlines                     \put(0,72){$c$}
       \put(19,57){$d$}
         \put(14,95){$c^{o}$}
          \put(5,72){\line(6,5){10}}
          \put(30,80){\line(1,0){10}}
          \put(20,80){$b$}
         \put(38,57){$e$}
         \put(38,83){$a$}
          \put(61,45){$e'$}
        \put(47,83){\line(6,5){10}}
        \put(57,93){\line(1,0){10}}
         \put(47,57){\line(6,-5){10}}
         \put(26,58){\line(1,0){10}}
     \put(7,72){\line(6,-5){10}}
         \put(7,75){\line(1,5){4}}
         \put(68,95){$(a')$} 
                \put(85,72){$z''_{3}$}
         \end{picture}\\ \cline{2 - 4}
& $z_{1}=(cc'bdaee')^{9}$ & $z_{2}=(cc'bdeaa')^{9}$ & $z_{3}=(ee'dfcaa')^{9}$\\
& $z'_{1}=(ee'dfacc')^{9}$ & $z'_{2}=(aa'bfecc')^{9}$ & $z'_{3}=(aa'bfcee')^{9}$\\
& $z''_{1}=(aa^{o}bfecc')^{9}$ & $z''_{2}=(ee^{o}dfacc')^{9}$ & $z''_{3}=(cc^{o}bdaee')^{9}$ \\
\hline\hline
 & $<m_{c},m_{e}>.<m_{c}z_{1}>$ & $<m_{c},m_{a}>.<m_{c}z_{2}>$ & $<m_{e},m_{a}>.<m_{e}z_{3}>$ \\ \cline{2 - 4}
  $Z(K)$ & $m_{c}m_{e}=n_{a}$ & $m_{c}m_{a}=n_{e}$ & $m_{a}m_{e}=n_{c}$ \\ \cline{2 - 4}
& $m_{c}z_{1}=m_{e}z'_{1}=n_{a}z''_{1}$ & $m_{c}z_{2}=m_{a}z'_{2}=n_{e}z''_{2}$ & $m_{e}z_{3}=m_{a}z'_{3}=n_{c}z''_{3}$ \\
\hline

\end{array}
\]
\pagebreak
\subsection{Lien entre $E$ et $G$}\label{13}
On note $\Phi$ la correspondance entre les générateurs de $G$ et des éléments de $C_{E}(d_{3})$ définie de la manière suivante:
\[
\begin{array}{  >{$}l<{$} >{$}r<{$}}
\begin{picture}(180,120)
\thinlines  
	\put(-20,72){$G:$}
          \put(8,72){$c'$}
          \put(30,72){$c$}
         \put(49,57){$d$}
        \put(49,83){$b$}
        \put(15,72){\line(1,0){12}}
        \put(68,57){$e$}
        \put(68,83){$a$}
        \put(88,72){$f$}
        \put(77,83){\line(6,5){10}}
       \put(77,83){\line(6,-5){10}}
        \put(77,61){\line(6,5){10}}
        \put(56,83){\line(1,0){10}}
        \put(56,58){\line(1,0){10}}
        \put(37,72){\line(6,-5){10}}
        \put(37,72){\line(6,5){10}}
       \put(77,57){\line(6,-5){10}}
        \put(88,45){$e'$}
        \put(88,90){$a'$}
        \end{picture}

 &\begin{picture}(140,120)
 \put(-6,72){$E:$}    
\thinlines  \put(12,72){$d_{3}$}
                      \put(38,72){$c_{3}$}
            \put(62,72){$b_{3}$}
            \put(106,60){$b_{1}$}
            \put(25,72){\line(1,0){10}}
          \put(156,50){$d_{1}$}
             \put(86,72){$a$}
     \put(118,60){\line(1,0){10}}
     \put(131,60){$c_{1}$}
       \put(142,84){\line(6,-5){10}}
                   
         \put(106,82){$b_{2}$} 
      \put(131,82){$c_{2}$}
           \put(118,82){\line(1,0){10}}
           \put(156,94){$d_{2}$}
       \put(142,62){\line(6,-5){10}}
      \put(142,62){\line(6,5){10}}
        \put(156,72){$t$}
      \put(142,84){\line(6,5){10}}
                    \put(50,72){\line(1,0){10}}
             \put(74,72){\line(1,0){10}}
       \put(94,72){\line(6,-5){10}}
              \put(94,72){\line(6,5){10}}
              \end{picture}

\end{array}
\]
\[
\begin{array}{>{$}c<{$} >{$}c<{$} >{$}c<{$}}
$a \to  c_{2}$ & $d \to b_{1}$ & $a' \to d_{2}$\\
$b \to b_{2}$ & $e \to c_{1}$ & $c' \to b_{3}$\\
$c \to a$ & $f \to t$ & $e' \to d_{1}$
\end{array}
\]

On note $D$ le sous-groupe de $E$ engendré par les éléments $\Phi (a)$, $\Phi (b)$, ....
\subsubsection{Hexagonale}\label{131}
La relation hexagonale est satisfaite par les images des générateurs de $G$. (\ref {H})
\subsubsection{Éléments centraux de $E$}\label{132}
Les éléments $\Phi (z_{1})$ et $\Phi(z_{2})$ écrits comme produits des éléments $\Phi(c) \ldots$ sont les éléments centraux $f_{1}$ et $f_{2}$ de $E$. On vérifie que $\Phi(z_{3})$ est central dans $E$. Ainsi la relation $\mathcal{R}$ est satisfaite dans $E$.

Les éléments $\Phi(m_{a})$ et $\Phi(m_{e})$ sont centraux dans $E$, ils n'appartiennent pas à l'ensemble
$Z(E) - \{1\}$: on a donc $\Phi(m_{a})=\Phi(m_{e})=1$.

Enfin on a $\Phi(m_{c})=m$ (notation \ref{E}); c'est un élément de $Z(D)$ qui n'est pas dans $Z(E)$.
\subsubsection{La relation $r=1$}\label{133}
Les éléments $\mathcal{C}^{ee'de}$ et $\mathcal{A}$ correspondent par $\Phi$ à $x=a^{b_{1}b_{2}b_{3}a.c_{1}d_{1}b_{1}c_{1}}$ et $y=c_{2}^{d_{2}b_{2}tc_{2}}$; la relation $r=1$, $(\mathcal{C}^{ee'de}. \mathcal{A})^{4}=1$, s'écrit donc dans $E$ $(xy)^{4}=1$. Rappelons que $E$ est un groupe de $\{3,4\}$-transpositions (\cite{[2]}, \cite{[6]}), l'ordre du produit de deux transpositions distinctes est donc $2,3,4$. Or on sait que $C_{E}(d_{3})/O_{2}(C_{E}(d_{3}))$ est isomorphe à $U_{6}(2)$, groupe dans lequel l'image de $xy$ n'est pas d'ordre $3$; on a $(xy)^{4}=1$, la relation $r=1$ est satisfaite dans $E$.
\subsubsection{En conclusion}\label{134}
De ce qui précède, il résulte que $\Phi$ induit un morphisme de $G$ dans $E$ dont l'image $D$ est un sous-groupe du centralisateur de $d_{3}$ dans $E$.
\section{Le sous-groupe $N$}
Dans cette section $G$ désigne un groupe satisfaisant aux hypothèses du théorème avec les relations $rel(3)$. Les notations sont celles introduites en \ref{not}.

L'objet de cette section est l'étude de la fermeture normale $N$ dans $G$ de $\alpha _{a'}=a'a^{o}$,  $\alpha _{a'}$ est un élément d'ordre 2 (voir \ref{zi}). On détermine un système générateur de $N$ de cardinal 22, on établit que le groupe des commutateurs de $N$ est d'ordre 2 et que le centre de $N$ est abélien élémentaire d'ordre 8 et que l'on a $\mathcal{D}(N)\subset Z(N)\subset Z(G)$. Enfin, on vérifie que $N$ est d'ordre $2^{23}$.

On pose $Y=\{a,b,\ldots,f,a',c'\}$; $Y\cup\{e'\}$ et $Y\cup\{a^{o}\}$ sont des systèmes générateurs de $G$.
\subsection{Les éléments $\alpha _{y}$ et $\beta _{y}$, $y\in Y$}
\subsubsection{} \label{211}
L'élément $\alpha_{a'}$ et certains de ses conjugués dans $G$ appartiennent au sous-groupe $W_{12}=<W_{1},W_{2}>$ (notations \ref{zi}); la fermeture normale de $\alpha_{a'}$ dans $W_{12}$ est un 2-groupe abélien élémentaire d'ordre $2^{7}$ contenant $z_{1}z_{2}$ (voir l'Annexe 3) et est engendré par les éléments $\alpha_{a'}$, $\alpha_{a}=\alpha_{a'}^{aa'}$, $\alpha_{b}=\alpha_{a}^{ba}$, $\alpha_{c}=\alpha_{b}^{cb}$, $\alpha_{c'}=\alpha_{c}^{c'c}$, $\alpha_{d}=\alpha_{c}^{dc}$, $\alpha_{e}=\alpha_{d}^{ed}$ et l'on a $\alpha_{a'}\alpha_{b}\alpha_{c'}=z_{1}z_{2}$ (voir aussi \cite{Z1}).

On pose $\alpha_{f}=\alpha_{e}^{fe}$ et $y_{1}=y\alpha_{y}$ pour $y$ dans $Y$; on vérifie que les produits $yy'$ et $y_{1}y'_{1}$ ont le même ordre pour $y$ et $y'$ dans $Y-\{f\}$. On calcule les ordres de $yf_{1}$ et $fy_{1}$ ($y$ dans $Y$), voir la table T.1.
\subsubsection{}\label{212}
On observe que $f_{1}^{af}=a_{1}^{bcv^{2}cb}$ avec $v=adbecf$; puisque $v^{4}=V=1$, on a
$a_{1}=a_{1}^{bcv^{2}cb.bcv^{2}cb}=f_{1}^{afbacbdcedfeaf}$. On pose $\beta_{a}=\alpha_{f}^{af}$, on définit successivement $\beta_{a'}=\beta_{a}^{a'a}$, $\beta_{b}=\beta_{a}^{ba}$, $\beta_{c}=\beta_{b}^{cb}$, $\beta_{c'}=\beta_{c}^{c'c}$, $\beta_{d}=\beta_{c}^{dc}$, $\beta_{e}=\beta_{d}^{ed}$, $\beta_{f}=\beta_{e}^{fe}$.
On a alors $\beta_{f}^{af}=\alpha_{a}$ d'où $\beta_{f}^{a}=\beta_{f}\alpha_{a}=\alpha_{a}\beta_{f}=\alpha_{a}^{f}$ et $\alpha_{f}^{a}=\alpha_{f}\beta_{a}=\beta_{a}\alpha_{f}=\beta_{a}^{f}$.
Pour $y$ dans $Y$, on pose $y_{2}=y\beta_{y}$.
\subsubsection{}\label{213}
On a ainsi déterminé 16 éléments de $N$; la table T.1 donne l'ordre des produits suivants: $yy'_{1}$, $yy'_{2}$, $y_{2}y'_{1}$, $y_{2}y'_{2}$ et $y_{1}y'_{1}$ pour $y$ et $y'$ dans $Y$ (où $yy_{1}$ et $yy_{2}$ désignent les involutions $\alpha_{y}$ et $\beta_{y}$).
\subsection{Les éléments $x_{a}$, $x_{a'}$, $x_{b}$, $x_{f}$}
\subsubsection{}\label{221}
Soit $B$ le sous-groupe de $G$ engendré par $a$, $a'$, $b$, $f$ et $a^{o}$. 
Posons $L_{a}=a^{a'bfa^{o}}a$, $L_{a'}=L_{a}^{a'a}$, $L_{b}=L_{a}^{ba}$, $L_{f}=L_{a}^{fa}$ et $L_{a^{o}}=L_{a}^{a^{o}a}$. Ces éléments sont d'ordre 4 (relation $r=1$); on désigne par $x_{a}$, $x_{a'}$, $x_{b}$, $x_{f}$ et $x_{a^{o}}$ leur carré respectif. Il est facile de voir que ces éléments commutent entre eux et que l'on a $x_{a}x_{a'}x_{b}x_{f}x_{a^{o}}=1$ (\cite{[8]}, \cite{Z1}, \cite{Z2}).

On observe que $x_{f}=(b_{1}f)^{2}=\alpha_{b}^{f}\alpha_{b}=\alpha_{b}\alpha_{b}^{f}$ et $x_{b}=(f_{2}b)^{2}=\beta_{f}^{b}\beta_{f}=\beta_{f}\beta_{f}^{b}$; ainsi les éléments $x_{a}$, $x_{a'}$, $x_{b}$, $x_{f}$ sont dans $N$ et l'on a $\alpha_{b}^{f}=\alpha_{b}x_{f}$ et $\beta_{f}^{b}=\beta_{f}x_{b}$.
\subsubsection{}\label{222}
Des égalités $\alpha_{a}=\alpha_{b}^{ab}=\beta_{f}^{af}$, $\beta_{c}=\beta_{b}^{cb}=\beta_{d}^{cd}$ et $\alpha_{e}=\alpha_{d}^{ed}=\alpha_{f}^{ef}$ il vient: $\beta_{f}^{b}=\alpha_{b}\beta_{f}\alpha_{b}^{f}$, $\beta_{d}^{b}=\beta_{b}\beta_{d}\beta_{b}^{d}$, $\alpha_{f}^{d}=\alpha_{d}\alpha_{f}\alpha_{d}^{f}$; de la première égalité on tire: $\beta_{f}\beta_{f}^{b}=\beta_{f}(\alpha_{b}\beta_{f})\alpha_{b}\alpha_{b}\alpha_{b}^{f}$ d'où $x_{b}=[\beta_{f},\alpha_{b}]x_{f}$. Ainsi $x_{f}x_{b}$ est un élément de $\mathcal{D}(N)$.
\subsubsection{}\label{223}
De l'égalité $\alpha_{a'}\alpha_{b}\alpha_{c'}=z_{1}z_{2}$ (\ref{211}) il résulte que $\alpha_{a'}\alpha_{b}\alpha_{c'}$ commute à $f$, on a donc $\alpha_{c'}^{f}=\alpha_{c'}x_{f}$ (car $x_{f}$ commute avec $\alpha_{a'}$ et $\alpha_{b}$). On en déduit en outre que $\beta_{f}$ commute à $\alpha_{b}$ et que $[\alpha_{a},\beta_{b}]$ est un élément de $Z(G)$. On a donc $x_{b}=x_{f}$. Désormais on pose $k:=[\alpha_{a},\beta_{b}]$, $k$ est un élément de $Z(G)$.
\subsubsection{}\label{224}
Les éléments $m_{a}$ et $m_{c}$ sont dans le centre de $G$ (\ref{127}), ils commutent avec $\alpha_{b}$. On en tire les égalités suivantes: $\alpha_{d}^{f}=\alpha_{d}x_{f}$, $\beta_{a'}^{d}=\beta_{a'}x_{a'}$,  $\beta_{a}^{d}=\beta_{a}x_{a}$, $\beta_{b}^{d}=\beta_{b}x_{f}$ ainsi que l'égalité $x_{a'}=x_{f}$.
\subsection{Quelques éléments de $\mathcal{D}(N)$; les conjugués de $x_{a},x_{f},\alpha_{y},\beta_{y}$ ($y\in Y$) par les éléments de $Y$}
\subsubsection{Les commutateurs de $\alpha_{y}$ et $\beta_{y}$ avec $x_{a},x_{f},\alpha_{y'},\beta_{y'}$ ($y,y'\in Y$)}\label{231}
On établit d'abord que $x_{a}$ et $x_{f}$  commutent avec $\alpha_{y}$ et $\beta_{y}$ ($y\in Y$), puis que pour $\{y,y'\}\ne \{a,f\}$ et $y,y'\in Y$ on a $[\alpha_{y},\alpha_{y'}]=[\beta_{y},\beta_{y'}]=1$ et que $[\alpha_{a},\alpha_{f}]=[\beta_{a},\beta_{f}]=k$ ($k$ introduit en \ref{223}). Enfin, pour $y$, $y'$ dans $Y$, on a $[\alpha_{y},\beta_{y'}]=k$ si $yy'$ est d'ordre 3 et $[\alpha_{y},\beta_{y'}]=1$ si $yy'$ est d'ordre 2.

Ces résultats s'obtiennent grâce aux relations $Q_{222}$ et à la table des produits T.1 (\ref{213}).
\subsubsection{}\label{232}
La table T.2 donne les expressions des conjugués de $\alpha_{y}$, $\beta_{y}$, $x_{a}$, $x_{f}$ ($y\in Y$) par les éléments de $Y$, expressions écrites comme produit de ces mêmes éléments $\alpha_{y}$, $\beta_{y}$, $x_{a}$, $x_{f}$. Presque tous ces conjugués peuvent être écrits de cette manière. Il reste $\alpha_{c}^{f}$, $\beta_{e}^{b}$ et des éléments notés $t_{i}$ ($1\leqslant i \leqslant 4$).

Comme $\alpha_{b}^{f}=x_{f}\alpha_{b}$ (\ref{221}), en conjuguant par $c$ puis par $e$ on obtient d'abord $\alpha_{b}^{f}\alpha_{c}^{f}=x_{f}^{c}\alpha_{b}\alpha_{c}$, d'où $x_{f}^{c}=\alpha_{c}^{f}x_{f}x_{c}$ ce qui donne $t_{1}$, puis $(\alpha_{b}^{f})^{e}=x_{f}^{e}\alpha_{b}$; or $\beta_{f}^{b}=\beta_{f}x_{f}=\beta_{f}\alpha_{b}\alpha_{b}^{f}$ (\ref{221}, \ref{231}), il vient $\beta_{f}^{be}=\beta_{f}^{b}.\beta_{e}^{b}=\beta_{f}\beta_{e}\alpha_{b}\alpha_{b}^{fe}$ d'où $x_{f}^{e}=\beta_{e}^{b}x_{f}\beta_{e}$ ce qui donne $t_{2}$.

Des égalités $x_{a}^{e}=x_{b}^{ab.e}=x_{f}^{e.ab}$ et $x_{a}^{c}=x_{f}^{af.c}=x_{f}^{c.af}$ on déduit, en utilisant $t_{2}$ et $t_{1}$ les expressions de $t_{4}=x_{a}^{e}=\beta_{e}^{b}\beta_{e}^{ba}x_{a}$ et $t_{3}=x_{a}^{c}=\alpha_{c}^{f}\alpha_{c}^{fa}x_{f}$.
\subsection{L'ensemble $\Gamma=\{\alpha_{y},\beta_{y},x_{a},x_{f},\alpha_{c}^{f},\alpha_{c}^{fa},\alpha_{c}^{fe},\alpha_{c}^{fae},\beta_{e}^{b},\beta_{e}^{ba}\\(y \in Y)\}$ engendre $N$; on a $\mathcal{D}(N)=<k>$}\label{24}
On construit deux tables: celle, notée T.3, des conjugués $\gamma^{y} (\gamma \in \Gamma, y \in Y)$ écrits en fonction des éléments de $\Gamma$, et celle, notée T.4, des commutateurs des éléments de $\Gamma$. Les éléments non connus de la table T.3 sont notés $t_{5},t_{6},\ldots,t_{20}$, leurs expressions sont données dans la liste T.5.
\subsubsection{}\label{241}
La table T.4 s'obtient dans presque tous les cas par des arguments simples (relations $Q_{2,2,2}$, tables T.1 et T.2, résultats connus de la table T.3). dans les autres cas, on montre que chaque commutateur  est central dans $N$ et dans $G$, qu'il est de carré 1; puis on vérifie qu'il est égal soit à 1 soit à $k$. En particulier on a $k^{2}=1$ et la relation $r=1$, $(\mathcal{C}^{ee'de}\mathcal{A})^{4}=1$, impose $k\ne 1$.
\subsubsection{}\label{242}
En conjuguant $x_{f}^{c}$ par $e$ et $x_{f}^{e}$ par $c$ on obtient $x_{f}^{ce}=(\alpha_{c}^{f}x_{f}\alpha_{c})^{e}=(\beta_{e}^{b}x_{f}\beta{e})^{c}$ (\ref{232}) d'où l'on déduit les expressions de $\beta_{e}^{bc}$ ($t_{17}$) et $\beta_{e}^{bca}$ ($t_{20}$).

En outre, on a $x_{a}^{ca'a}=x_{f}^{c}$ puisque $x_{a'}=x_{f}$ (\ref{224}); on en tire les valeurs de $\alpha_{c}^{faa'}$ ($t_{5}$) et de $\alpha_{c}^{faa'e}$ ($t_{12}$). En écrivant $x_{a}^{ea'}=x_{a}^{e}x_{f}^{e}$, il vient l'expression de $\beta_{e}^{baa'}$ ($t_{19}$).
\subsubsection{}\label{243}
En utilisant les valeurs connues de la table T.2 on obtient par conjugaison les expressions de $\alpha_{c}^{fac'}$ ($t_{7}$), $\alpha_{c}^{fad}$ ($t_{8}$), $\alpha_{c}^{feb}$ ($t_{9}$), $\alpha_{c}^{fec'}$ ($t_{10}$), $\alpha_{c}^{faec'}$ ($t_{14}$), $\beta_{e}^{bd}$ ($t_{18}$), $\beta_{e}^{bad}$ ($t_{21}$).
\subsubsection{}\label{244}
Les éléments $m_{c}$ et $m_{a}$ sont dans $Z(G)$ (\ref{127}). De $[m_{c},\alpha_{e}]=1$ on déduit les écritures de $\alpha_{c}^{fab}$ ($t_{6}$) puis de $\alpha_{c}^{fabe}$ ($t_{13}$) et de $[m_{c},\alpha_{a}]=1$ on déduit celles de $\alpha_{c}^{fed}$ ($t_{11}$) et de $\alpha_{c}^{faed}$ ($t_{15}$) .

La connaissance de $t_{6}$, $t_{15}$ et $t_{20}$ permet d'écrire le conjugué de $\beta_{e}^{ba}$ par $\mathcal{C}=c^{c'bdc}$: 
\[
\beta_{e}^{ba\mathcal{C}}=kx_{a}x_{f}\alpha_{c'}\beta_{b}\beta_{f}\beta_{e}^{ba}. 
\]
L'égalité $[m_{c},\beta_{a}]=1$ conduit à une relation (S1) entre les éléments de $\Gamma$:\\
\begin{center}
(S1): \qquad $\alpha_{b}\alpha_{d}\alpha_{f}\beta_{b}\beta_{d}\beta_{f}=1$
\end{center}
Les expressions de $t_{6}$, $t_{12}$ et $t_{13}$ conduisent à une expression du conjugué de $\alpha_{c}^{fe}$ par $\mathcal{A}=a^{a'bfa}$:
\[
\alpha_{c}^{fe.\mathcal{A}}=k\alpha_{a'}\alpha_{d}\beta_{a'}\beta_{d}\beta_{f}\beta_{e}^{ba}x_{a}\alpha_{c}^{fe}\beta_{e}^{baf}.
\]

La relation (S1), les égalités $[m_{a},\beta_{e}]=[m_{a},\alpha_{c}]=1$ conduisent à l'expression de $\beta_{e}^{baf}$ ($t_{22}$). Par conjugaison de $t_{22}$ par $c$ et de $t_{20}$ par $f$ on obtient $\beta_{e}^{baf.c}=\beta_{e}^{bac.f}$ d'où:
\[
(\beta_{c'}\beta_{f}\beta_{e}x_{a}x_{f})^{c}\beta_{e}^{bac}\beta_{e}^{bc}=\alpha_{c}^{faef}\beta_{e}^{baf}\alpha_{c}^{faf}.
\]
Grâce aux expressions $t_{20}$, $t_{17}$ et $t_{22}$ on en déduit
\[
\alpha_{c}^{faef}=k\beta_{c}\alpha_{c}\alpha_{c}^{fae}\alpha_{c}^{fe}\alpha_{c}^{f}\alpha_{c}^{fa} \quad (t_{16})
\]
Enfin la relation (S1) et l'égalité $[m_{a},\alpha_{e}]=1$ conduisent à une seconde relation entre les éléments de $\Gamma$:\\
\begin{center}
(S2):\qquad $\alpha_{a'}\beta_{c'}\beta_{d}\beta_{f}=1$.
\end{center}
Les tables T.3 et T.4 sont achevées.
\subsubsection{}\label{245}
Le sous-groupe $N$ est engendré par $\Gamma$ (tables T.2, T.3) et son groupe des commutateurs $\mathcal{D}(N)$ est engendré par l'élément $k$ d'ordre 2 (table T.4). Les relations (S1) et (S2) prouvent qu'il y a dans $\Gamma$ deux générateurs superflus, par exemple $\beta_{b}=\alpha_{b}\alpha_{d}\alpha_{f}\beta_{d}\beta_{f}$ \quad (S1) et $\beta_{c'}=\alpha_{a'}\beta_{d}\beta_{f}$ \quad (S2). Ainsi $N$ est engendré par 22 involutions, l'ordre de $N/\mathcal{D}(N)$ divise $2^{22}$ et celui de $N$ divise  $2^{23}$.
\subsubsection{}\label{246}
On pose $\Gamma_{0} =\Gamma - \{\beta_{b},\beta_{c'}\}$. Toute relation entre les éléments de $\Gamma_{0}$ s'écrit $1=k^{p_{o}}\prod_{\gamma \in \Gamma_{o}}\gamma^{p_{\gamma}}$ où $p_{0}$ et $p_{\gamma}$ sont dans $\{0,1\}$, chaque élément de $\Gamma_{0}$ intervenant au plus une fois. On obtient alors: $(\alpha_{a'}\alpha_{b}\alpha_{c'})^{p}=1$ et $(\alpha_{c'}\alpha_{d}\alpha_{f}\beta_{a'})^{q}=1$ avec $p$ et $q$ dans $\{0,1\}$.
\subsubsection{}\label{247}
Les éléments $m_{c}z_{1}=m_{e}z'_{1}$, $m_{c}z_{2}=m_{a}z'_{2}$ et $m_{e}z_{3}=m_{a}z'_{3}$ sont respectivement des éléments centraux des sous-groupes $K_{a}$, $K_{e}$ et $K_{c}$ (voir \ref{126}). Des calculs conduisent aux égalités:
\begin{center}
\[
z_{1}z_{2}=\alpha_{a'}\alpha_{b}\alpha_{c'}, \quad z'_{1}z_{3}=\alpha_{c'}\alpha_{d} \alpha_{f}\beta_{a'},\quad z'_{2}z'_{3}=\alpha_{a'}\alpha_{b}\alpha_{d}\alpha_{f}\beta_{a'}.
\]   
\end{center}
On observe que les éléments $z_{1}z_{2}$ et $z'_{1}z_{3}=m_{e}m_{c}z_{1}.z_{3}$ sont dans $N\bigcap Z(G)$; on pose $z_{1}z_{2}=z$ et $z'_{1}z_{3}=\hat{z}$, on a $z'_{2}z'_{3}=z\hat{z}$.
\subsection{Le centre de $N$, l'ordre de $N$}\label{25}
\subsubsection{}\label{251}
Les éléments $m_{a}$, $m_{c}$ et $m_{e}$ sont centraux dans $G$ et l'on a $m_{a}m_{c}=n_{e}$, $m_{c}m_{e}=n_{a}$, $m_{e}m_{a}=n_{c}$ (\ref{126}).
Rappelons que $m_{a}=a'bfd\mathcal{A}\mathcal{A}^{cbdc}$ et $n_{a}=a'bfd\mathcal{A}_{0}\mathcal{A}_{0}^{cbdc}$ où $\mathcal{A}$ et $\mathcal{A}_{0}$ désignent respectivement $a^{a'bfa}$ et $a^{a^{o}bfa}$ (voir \ref{Ke} et \ref{126}). On établit alors l'égalité $m_{a}n_{a}=k$ et par suite $m_{c}n_{c}=k=m_{e}n_{e}$.
\subsubsection{}\label{252}
Il résulte de \ref{246} et \ref{247} que les éléments $k$, $z$ et $\hat{z}$ engendrent le centre de $N$ et que $Z(N)$ est contenu dans $Z(G)$. Ainsi on a $|Z(N)|=2^{3}$ et $|N|=2^{23}$.
\section{La preuve du théorème}
Pour $1\leqslant i \leqslant 3$, on considère un groupe $G_{i}$ avec la présentation
\[ 
(a,b,\ldots,f,a', c',e'/Q_{222},V=1,rel(i)) 
\]
(notations \ref{01}); on note $N_{i}$ la fermeture normale de $a'a^{o}$ dans $G_{i}$ (pour $a^{o}$ voir \ref{zi}) et $H_{i}$ le sous-groupe de $G_{i}$ engendré par $a,b,\ldots,f,a',c'$. Pour chaque jeu de relations $ rel(i)$, on précise la structure de $N_{i}$ (\ref{31}), celle de $H_{i}$ (\ref{32}) puis on établit que $G_{i}/Z(G_{i})N_{i}$ est isomorphe à $U_{6}(2)$ (\ref{33}). Enfin on fait le lien avec le centralisateur d'une involution de la classe $2A$ de $^{2}E_{6}(2)$.
\subsection{Le sous-groupe $N_{i}$ de $G_{i}$ }\label{31}
On a établi que $k$ (\ref{223}) est une involution centrale de $G_{i}$ qui engendre $\mathcal{D}(N_{i})$ et appartient à $Z(N_{i})$ (\ref{223}, \ref{241}). D'après \ref{247} et \ref{25}, le centre de $N_{i}$ est engendré par $k$, $z=z_{1}z_{2}$ et $\hat{z}=z'_{1}z_{3}$ (avec $z'_{1}=m_{e}m_{c}z_{1}$).

Sous $rel(1)$, on a $z_{1}=z_{2}=z_{3}=1$ et $m_{a}=m_{e}=1$; comme $m_{a}m_{e}=n_{c}=km_{c}$ (\ref{251}), il vient $k=m_{c}=\hat{z}$ et $z=1$. Ainsi on a $\mathcal{D}(N_{1})=Z(N_{1})=<k>$, $N_{1}$ est un $2$-groupe extraspécial. De plus l'ensemble générateur $\Gamma_{0}$ (\ref{246}) comporte deux éléments superflus: $\alpha_{a'}=\alpha_{b}\alpha_{c'}$ ($z=1$) et $\beta_{a'}=\alpha_{c'}\alpha_{d}\alpha_{f}k$ ($\hat{z}k=1$). Le sous-groupe $N_{1}$ est engendré par 20 éléments, $N_{1}$ est d'ordre $2^{20+1}$, c'est le produit de 10 groupes diedraux d'ordre 8 de même centre $<k>$:\\$<\alpha_{a},\alpha_{f}>.<\alpha_{b},\beta_{a}>.<\alpha_{c'},\beta_{c}>.<\alpha_{c},\beta_{d}>.<\alpha_{d},\beta_{e}>.<\alpha_{e},\beta_{f}>.\\<\alpha_{c}^{fa},\beta_{e}^{d}>.<\alpha_{c}^{f},\beta_{e}^{ba}>.<x_{f},\alpha_{c}^{fae}\alpha_{c}^{f}>.<x_{a},\alpha_{c}^{fe}\alpha_{c}^{f}>$.

Sous $rel(2)$ et $rel(3)$,  $N_{i}$  d'ordre $2^{23}$ (\ref{252}).
\subsection{Le sous-groupe $H_{i}$ de $G_{i}$}\label{32}
C'est un sous-groupe isomorphe à un quotient de $2^{2}.2.U_{6}(2)$ (voir Annexe 4) dont le centre est $<m_{a},m_{c}>.<m_{c}z_{2}>$ (Annexe 4 et \ref{Ke}); de plus $m_{a}m_{c}=km_{e}$ (\ref{126} et \ref{251}).

Sous $rel(1)$, on a $z_{2}=m_{a}=m_{e}=1$ en conséquence $m_{c}=k$ est l'unique involution centrale de $H_{1}$: $H_{1}$ est isomorphe à $2.U_{6}(2)$.

Sous $rel(2)$, $m_{c}=k$ et $z_{2}$ sont des involutions centrales de $H_{2}$; $Z(H_{2})$ est d'ordre 4 et $H_{2}$ est isomorphe à $2.2.U_{6}(2)$.

Sous $rel(3)$, $H_{3}$ est isomorphe à $2^{2}.2.U_{6}(2)$.
\subsection{Le groupe $G_{i}$}\label{33}
Soit $T_{i}$ le sous-groupe de $Z(G_{i})$ engendré par $z_{1}$,  $z_{2}$,  $z_{3}$, $m_{a}$, $m_{c}$, $m_{e}$, et soit $\pi_{i}$ l'application canonique $G_{i} \to G_{i}/T_{i}N_{i}$. Remarquons que $m_{a}m_{e}=km_{c}$, $z=z_{1}z_{2}$ et $\hat{z}=z'_{1}z_{3}=m_{e}m_{c}z_{1}z_{3}$ sont dans $T_{i}$, donc $Z(N_{i})$ est un sous-groupe de $T_{i}$ dont l'indice est 1, 2 ou $2^{3}$ suivant que $i=1$, $i=2$ ou $i=3$.

L'image de $H_{i}$ par $\pi_{i}$ est un sous-groupe $U_{i}$ isomorphe à $U_{6}(2)$ ($Z(H_{i})\subset T_{i}$) et l'image du système générateur de $G_{i}$ coïncide avec celle de $H_{i}$ (on a $\pi_{i}(a')=\pi_{i}(a^{o}$)). Par conséquent $\pi_{i}(G_{i})$ est un quotient de $2^{2}.2.U_{6}(2)$. Mais comme $\pi_{i}(G_{i})=\pi_{i}(H_{i})=U_{i}$, dans toutes les situations $i=1$, $i=2$, $i=3$, $\pi_{i}(G_{i})$ est isomorphe à $U_{6}(2)$ et l'on a $Z(G_{i})=T_{i}$. On en déduit les isomorphismes suivants:
\begin{center}
$G_{1}\simeq 2^{20+1}.U_{6}(2), \quad G_{2}\simeq 2^{24}.U_{6}(2), \quad G_{3}\simeq 2^{2}.2^{24}.U_{6}(2)$.
\end{center}
\subsection{Fin de la preuve}\label{34}
Rappelons que le groupe $E=2^{3}.^{2}E_{6}(2)$ admet la présentation
\begin{center}
 $(a,b_{i},c_{i},d_{i}\quad (1\leqslant i \leqslant 3)/Q_{222},V=1)$:
\end{center}
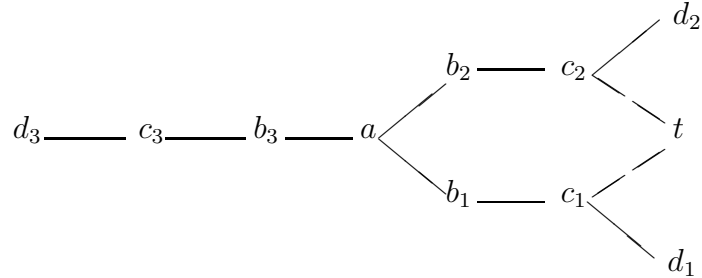
\[
\begin{picture}(150,150)
\thinlines  \put(-30,72){$d_{3}$}
                      \put(17,72){$c_{3}$}
            \put(60,72){$b_{3}$}
            \put(132,48){$b_{1}$}
            \put(-18,72){\line(1,0){30}}
          \put(215,22){$d_{1}$}
             \put(100,72){$a$}
     \put(144,48){\line(1,0){25}}
     \put(175,48){$c_{1}$}
      \put(185,48){\line(6,-5){25}}
                   
         \put(132,96){$b_{2}$} 
      \put(175,96){$c_{2}$}
           \put(144,98){\line(1,0){25}}
           \put(217,115){$d_{2}$}
        \put(187,96){\line(3,-2){12}}
       \put(202,86){\line(3,-2){12}}
      \put(187,50){\line(3,2){12}}
      \put(202,60){\line(3,2){12}}
        \put(217,72){$t$}
    \put(187,96){\line(6,5){25}}
                    \put(27,72){\line(1,0){30}}
             \put(72,72){\line(1,0){25}}
       \put(107,72){\line(6,-5){25}}
              \put(107,72){\line(6,5){25}}
              \end{picture}
\]
\begin{center}
$V=1$ \qquad avec $V=(atb_{1}c_{2}c_{1}b_{2})^{4}$.
\end{center}
Sous $rel(2)$ (resp. $rel(3)$) on a un morphisme $\Phi_{i}$ de $G_{i}$ dans le centralisateur de $d_{3}$ dans $E$ qui envoie $z_{1},z_{2},z_{3}$ dans $Z(E)$ (voir \ref{13});
 l'image de $G_{2}$ (resp. $G_{3}$) est un sous-groupe du centralisateur de $d_{3}$ dans $E$ et l'on a $|C_{E}(d_{3})|=|G_{2}|$. Ainsi \[(a,b,\ldots,f,a',c',e'/Q_{222},V=1, rel(i))\] est une présentation de $C_{E}(d_{3})$ pour $i=2$ et d'une extension $2^{2}.C_{E}(d_{3})$ pour $i=3$.

Sous $rel(1)$, on a un morphisme de $G_{1}$ dans $\bar{E}=E/Z(E)$ qui envoie $G_{1}$ dans le centralisateur $C$ de l'image de $d_{3}$ dans $\bar{E}$; comme $G_{1}$ et $\mathcal{C}$ ont le même ordre, 
\begin{center}
$(a,b,\ldots,f,a',c',e'/Q_{222},V=1,rel(1))$
\end{center}
est une présentation du centralisateur d'une involution $d$ dans 
$^{2}E_{6}(2)$ ($d$ provenant de la classe $2A$).
\section{Tables.}
Table des ordres des produits $yy'_{i}$, $y_{2}y'_{i}$, $y_{i}y'_{i}$ pour $i=1,2$ et $y,y'\in Y$.
\[
\begin{array}{r|cccccccc||cccccccc||}
&a'_{1}&a_{1}&b_{1}&c_{1}&c'_{1}&d_{1}&e_{1}&f_{1}&a'_{2}&a_{2}&b_{2}&c_{2}&c'_{2}&d_{2}&e_{2}&f_{2} \\ \hline
a'&2&3&2&2&2&2&2&2&2&3&2&2&2&2&2&2 \\ \hline
a&3&2&3&2&2&2&2&3&3&2&3&2&2&2&2&3 \\ \hline
b&2&3&2&3&2&2&2&2&2&3&2&3&2&4&4&4 \\ \hline
c&2&2&3&2&3&3&2&2&2&2&3&2&3&3&2&2 \\ \hline
c'&2&2&2&3&2&2&2&2&2&2&2&3&2&2&2&2 \\ \hline
d&2&2&2&3&2&2&3&4&4&4&4&3&2&2&3&2 \\ \hline
e&2&2&2&2&2&3&2&3&2&2&2&2&2&3&2&3 \\ \hline
f&2&3&4&4&4&4&3&2&2&3&2&2&2&2&3&2 \\ \hline \hline
a'_{2}&2&3&2&2&2&4&2&2&2&3&2&2&2&4&2&2 \\ \hline
a_{2}&3&2&3&2&2&4&2&3&&2&3&2&2&4&2&3 \\ \hline
b_{2}&2&3&2&3&2&4&2&2&&&2&3&2&2&4&4 \\ \hline
c_{2}&2&2&3&2&3&3&2&2&&&&2&3&3&2&2 \\ \hline
c'_{2}&2&2&2&3&2&2&2&2&&&&&2&2&2&2 \\ \hline
d_{2}&2&2&4&3&2&2&3&4&&&&&&2&3&2 \\ \hline
e_{2}&2&2&4&2&2&3&2&3&&&&&&&2&3 \\ \hline
f_{2}&2&3&2&4&4&4&3&2&&&&&&&&2 \\ \hline \hline
a'_{1}&2&3&2&2&2&2&2&2&&&&&&&& \\ \hline
a_{1}&&2&3&2&2&2&2&3&&&&&&&& \\ \hline
b_{1}&&&2&3&2&2&2&4&&&&&&&& \\ \hline
c_{1}&&&&2&3&3&2&4&&&&&&&& \\ \hline
c'_{1}&&&&&2&2&2&4&&&&&&&& \\ \hline
d_{1}&&&&&&2&3&2&&&&&&&& \\ \hline
e_{1}&&&&&&&2&3&&&&&&&& \\ \hline
f_{1}&&&&&&&&2&&&&&&&& \\ \hline \hline
\end{array}
\]
\begin{center}
Table T.1
\end{center}
\newpage

\begin{center}
Tables des conjugués  $\gamma_{y}^{y'}$ pour $\gamma_{y}\in \{\alpha_{y},\beta_{y}\}$, $y$ et $y'\in Y$
\end{center}
\begin{center}
Table T.2 (\ref{232})
\end{center}
\[
\begin{array}{|c |c |c |c |c |c |c |c |c| }\hline
&a'&a&b&c&c'&d&e&f \\ \hline
\alpha_{a'}&\alpha_{a'}&\alpha_{a'}\alpha_{a}&\alpha_{a'}&\alpha_{a'}&\alpha_{a'}&\alpha_{a'}&\alpha_{a'}&\alpha_{a'}\\ 
\alpha_{a}&\alpha_{a}\alpha_{a'}&\alpha_{a}&\alpha_{a}\alpha_{b}&\alpha_{a}&\alpha_{a}&\alpha_{a}&\alpha_{a}&\alpha_{a}\beta_{f} \\
\alpha_{b}&\alpha_{b}&\alpha_{b}\alpha_{a}&\alpha_{b}&\alpha_{b}\alpha_{c}&\alpha_{b}&\alpha_{b}&\alpha_{b}&\alpha_{b}x_{f} \\
\alpha_{c}&\alpha_{c}&\alpha_{c}&\alpha_{c}\alpha_{b}&\alpha_{c}&\alpha_{c}\alpha_{c'}&\alpha_{c}\alpha_{d}&\alpha_{c}&\alpha_{c}^{f} \\
\alpha_{c'}&\alpha_{c'}&\alpha_{c'}&\alpha_{c'}&\alpha_{c'}\alpha_{c}&\alpha_{c'}&\alpha_{c'}&\alpha_{c'}&\alpha_{c'}x_{f} \\
\alpha_{d}&\alpha_{d}&\alpha_{d}&\alpha_{d}&\alpha_{d}\alpha_{c}&\alpha_{d}&\alpha_{d}&\alpha_{d}\alpha_{e}&\alpha_{d}x_{f} \\
\alpha_{e}&\alpha_{e}&\alpha_{e}&\alpha_{e}&\alpha_{e}&\alpha_{e}&\alpha_{e}\alpha_{d}&\alpha_{e}&\alpha_{e}\alpha_{f} \\
\alpha_{f}&\alpha_{f}&\alpha_{f}\beta_{a}&\alpha_{f}&\alpha_{f}&\alpha_{f}&\alpha_{f}x_{f}&\alpha_{f}\alpha_{e}&\alpha_{f} \\ \hline
\beta_{a'}&\beta_{a'}&\beta_{a'}\beta_{a}&\beta_{a'}&\beta_{a'}&\beta_{a'}&\beta_{a'}x_{f}&\beta_{a'}&\beta_{a'} \\
\beta_{a}&\beta_{a}\beta_{a'}&\beta_{a}&\beta_{a}\beta_{b}&\beta_{a}&\beta_{a}&\beta_{a}x_{a}&\beta_{a}&\beta_{a}\alpha_{f} \\
\beta_{b}&\beta_{b}&\beta_{b}\beta_{a}&\beta_{b}&\beta_{b}\beta_{c}&\beta_{b}&\beta_{b}x_{f}&\beta_{b}&\beta_{b} \\
\beta_{c}&\beta_{c}&\beta_{c}&\beta_{c}\beta_{b}&\beta_{c}&\beta_{c}\beta_{c'}&\beta_{c}\beta_{d}&\beta_{c}&\beta_{c} \\
\beta_{c'}&\beta_{c'}&\beta_{c'}&\beta_{c'}&\beta_{c'}\beta_{c}&\beta_{c'}&\beta_{c'}&\beta_{c'}&\beta_{c'} \\
\beta_{d}&\beta_{d}&\beta_{d}&\beta_{d}x_{f}&\beta_{d}\beta_{c}&\beta_{d}&\beta_{d}&\beta_{d}\beta_{e}&\beta_{d} \\
\beta_{e}&\beta_{e}&\beta_{e}&\beta_{e}^{b}&\beta_{e}&\beta_{e}&\beta_{e}\beta_{d}&\beta_{e}&\beta_{e}\beta_{f} \\
\beta_{f}&\beta_{f}&\beta_{f}\alpha_{a}&\beta_{f}x_{f}&\beta_{f}&\beta_{f}&\beta_{f}&\beta_{f}\beta_{e}&\beta_{f} \\ \hline
x_{f}&x_{f}&x_{f}x_{a}&x_{f}&t_{1}&x_{f}&x_{f}&t_{2}&x_{f} \\
x_{a}&x_{a}x_{f}&x_{a}&x_{a}x_{f}&t_{3}&x_{a}&x_{a}&t_{4}&x_{a}x_{f} \\ \hline \hline\hline
\alpha_{c}^{f}&\alpha_{c}^{f}&\alpha_{c}^{fa}&\alpha_{c}^{f}\alpha_{b}x_{f}&\alpha_{c}^{f}&\alpha_{c}^{f}\alpha_{c'}x_{f}&\alpha_{c}^{f}\alpha_{d}x_{f}&\alpha_{c}^{fe}&\alpha_{c} \\
\alpha_{c}^{fa}&t_{5}&\alpha_{c}^{f}&t_{6}&\alpha_{c}^{fa}&t_{7}&t_{8}&\alpha_{c}^{fae}&\alpha_{c}^{fa}\\
\alpha_{c}^{fe}&\alpha_{c}^{fe}&\alpha_{c}^{fae}&t_{9}&\alpha_{c}^{fe}&t_{10}&t_{11}&\alpha_{c}^{f}&\alpha_{c}^{fe} \\
\alpha_{c}^{fae}&t_{12}&\alpha_{c}^{fe}&t_{13}&\alpha_{c}^{fae}&t_{14}&t_{15}&\alpha_{c}^{fa}&t_{16} \\
\beta_{e}^{b}&\beta_{e}^{b}&\beta_{e}^{ba}&\beta_{e}&t_{17}&\beta_{e}^{b}&t_{18}&\beta_{e}^{b}&\beta_{e}^{b}\beta_{f}x_{f} \\
\beta_{e}^{ba}&t_{19}&\beta_{e}^{b}&\beta_{e}^{ba}&t_{20}&\beta_{e}^{ba}&t_{21}&\beta_{e}^{ba}&t_{22} \\ \hline
\end{array}
\]
\begin{center}
Table T.3 (\ref{24})
\end{center}
\newpage
\begin{center}
Table des commutateurs des éléments de $\Gamma$ (\ref{24})\\
Table T.4
\end{center}
\[
\begin{array}{|c| c  c  c  c  c  c  c  c|c|c  c  c  c |c  c|}\hline
&\beta_{a'}&\beta_{a}&\beta_{b}&\beta_{c}&\beta_{c'}&\beta_{d}&\beta_{e}&\beta_{f}&\alpha_{f}&\alpha_{c}^{f}&\alpha_{c}^{fa}&\alpha_{c}^{fe}&\alpha_{c}^{fae}&\beta_{e}^{b}&\beta_{e}^{ba} \\ \hline
\alpha_{a'}&&k&&&&&&&&&&&&&\\ 
\alpha_{a}&k&&k&&&&&&k&&&&&&\\
\alpha_{b}&&k&&k&&&&&&&&&&&\\ \hline
\alpha_{c}&&&k&&k&k&&&&&&&&& \\
\alpha_{c'}&&&&k&&&&&&&&&&&\\ \hline
\alpha_{d}&&&&k&&&k&&&&&&&k&k \\
\alpha_{e}&&&&&&k&&k&&&&&&& \\
\alpha_{f}&&&&&&&k&&&&&&&k&k \\ \hline
\beta_{a'}&&&&&&&&&&&&&&& \\
\beta_{a}&&&&&&&&k&&&&&&& \\
\beta_{b}&&&&&&&&&&k&k&k&k&& \\ \hline
\beta_{c}&&&&&&&&&&&&&&& \\
\beta_{c'}&&&&&&&&&&k&k&k&k&& \\ \hline
\beta_{d}&&&&&&&&&&k&k&k&k&& \\ 
\beta_{e}&&&&&&&&&k&&&&&& \\
\beta_{f}&&k&&&&&&&&&&&&& \\ \hline
x_{a}&&&&&&&&&&&&k&k&& \\
x_{f}&&&&&&&&&&&&&k&& \\ \hline
\alpha_{c}^{f}&&&k&&k&k&&&&&&&&&k \\
\alpha_{c}^{fa}&&&k&&k&k&&&&&&&&k& \\ 
\alpha_{c}^{fe}&&&k&&k&k&&&&&&&&&k \\
\alpha_{c}^{fae}&&&k&&k&k&&&&&&&&k& \\ \hline
\beta_{e}^{b}&&&&&&&&&&&k&&k&& \\
\beta_{e}^{ba}&&&&&&&&&&k&&k&&& \\ \hline
\end{array}
\]
\begin{center}
Les valeurs non indiquées valent 1.
\end{center}
\newpage
\begin{center}
Valeurs des éléments $t_{i}$ ($1\leqslant i \leqslant 22$)\\ Table T.5
\end{center}
\[
\begin{array}{|c| c|c|c|}\hline
t_{1}&x_{f}^{c}&\alpha_{c}\alpha_{c}^{f}x_{f}&\ref{232} \\ \hline
t_{2}&x_{f}^{e}&\beta_{e}\beta_{e}^{b}x_{f}&\ref{232} \\ \hline
t_{3}&x_{a}^{c}&\alpha_{c}^{f}\alpha_{c}^{fa}x_{a}&\ref{232} \\ \hline
t_{4}&x_{a}^{e}&\beta_{e}^{b}\beta_{e}^{ba}x_{a}&\ref{232} \\ \hline
t_{5}&\alpha_{c}^{faa'}&\alpha_{c}\alpha_{c}^{f}\alpha_{c}^{fa}&\ref{242} \\ \hline
t_{6}&\alpha_{c}^{fab}&k\alpha_{c}\alpha_{d}\beta_{c'}\beta_{d}\beta_{f}\alpha_{c}^{f}\alpha_{c}^{fa}x_{a}x_{f}&\ref{244} \\ \hline
t_{7}&\alpha_{c}^{fac'}&\alpha_{c'}\alpha_{c}^{fa}x_{a}x_{f}&\ref{243} \\ \hline
t_{8}&\alpha_{c}^{fad}&\alpha_{d}\alpha_{c}^{fa}x_{a}x_{f}&\ref{243} \\ \hline
t_{9}&\alpha_{c}^{feb}&\alpha_{b}\beta_{e}\alpha_{c}^{fe}\beta_{e}^{b}x_{f}&\ref{243} \\ \hline
t_{10}&\alpha_{c}^{fec'}&\alpha_{c'}\beta_{e}\alpha_{c}^{fe}\beta_{e}^{b}x_{f}&\ref{243} \\ \hline
t_{11}&\alpha_{c}^{fed}&\alpha_{b}\alpha_{c}\alpha_{f}\beta_{b}\beta_{c'}\beta_{e}\alpha_{c}^{f}\alpha_{c}^{fe}\beta_{e}^{b}x_{f}&\ref{244} \\ \hline
t_{12}&\alpha_{c}^{faea'}&\alpha_{c}\alpha_{c}^{fe}\alpha_{c}^{fae}&\ref{242} \\ \hline
t_{13}&\alpha_{c}^{faeb}&k\alpha_{c}\alpha_{d}\alpha_{e}\beta_{c'}\beta_{d}\beta_{e}\beta_{f}\alpha_{c}^{fe}\alpha_{c}^{fae}\beta_{e}^{ba}x_{a}x_{f}&\ref{244} \\ \hline
t_{14}&\alpha_{c}^{faec'}&\alpha_{c'}\beta_{e}\alpha_{c}^{fae}\beta_{e}^{ba}x_{a}x_{f}&\ref{243} \\ \hline
t_{15}&\alpha_{c}^{faed}&\alpha_{a}\alpha_{b}\alpha_{c}\alpha_{f}\beta_{b}\beta_{c'}\beta_{e}\alpha_{c}^{fa}\alpha_{c}^{fae}\beta_{e}^{ba}x_{a}x_{f}&\ref{244} \\ \hline
t_{16}&\alpha_{c}^{faef}&k\alpha_{c}\beta_{c}\alpha_{c}^{f}\alpha_{c}^{fa}\alpha_{c}^{fe}\alpha_{c}^{fae}&\ref{244} \\\hline
t_{17}&\beta_{e}^{bc}&\alpha_{c}^{f}\alpha_{c}^{fe}\beta_{e}^{b}&\ref{242} \\\hline
t_{18}&\beta_{e}^{bd}&\beta_{d}\beta_{e}^{b}x_{f}&\ref{243} \\\hline
t_{19}&\beta_{e}^{baa'}&\beta_{e}\beta_{e}^{b}\beta_{e}^{ba}&\ref{242} \\\hline
t_{20}&\beta_{e}^{bac}&\alpha_{c}^{fa}\alpha_{c}^{fae}\beta_{e}^{ba}&\ref{242} \\\hline
t_{21}&\beta_{e}^{bad}&\beta_{d}\beta_{e}^{ba}x_{a}x_{f}&\ref{243} \\\hline
t_{22}&\beta_{e}^{baf}&\beta_{c'}\beta_{e}\beta_{f}\beta_{e}^{b}\beta_{e}^{ba}x_{a}x_{f}&\ref{244}\\\hline
\end{array}
\]
\section{Annexe}
Tous les graphes ci-dessous sont des graphes de Coxeter.
\subsection{Le graphe $Q_{111}$}
On désigne par $G$ le groupe $H_{3,6}\simeq 3^{5}\rtimes S_{6}$.
\subsubsection{}
Le groupe $G$ admet la présentation $(x_{i} (0\leqslant i \leqslant 5)/h)$ avec $h$:
\[
\begin{picture}(150,100)

\thinlines  \put(-30,72){$x_{1}$}
                      \put(17,72){$x_{2}$}
            \put(60,72){$x_{3}$}
           \put(-22,72){\line(1,0){30}}
          
             \put(100,72){$x_{4}$}
          \put(110,72){\line(1,0){23}}
    
       \put(-25,68){\line(5,-4){20}}
                
      \put(-2,48){$x_{0}$} 
              \put(10,52){\line(5,6){13}}
        \put(135,72){$x_{5}$}
   
                    \put(27,72){\line(1,0){30}}
              \put(72,72){\line(1,-0){25}}
              \put(-30,25){$(x_{1}^{x_{2}}x_{0})^{3}=1$}
              \end{picture}
\]
L'élément $z=(x_{1}x_{0})(x_{1}x_{0})^{x_{2}}(x_{1}x_{0})^{x_{2}x_{3}}(x_{1}x_{0})^{x_{2}x_{3}x_{4}}(x_{1}x_{0})^{x_{2}x_{3}x_{4}x_{5}}$ qui s'écrit aussi $(x_{2}x_{1}x_{3}x_{4}x_{5}x_{0})^{5}$ est d'ordre 3 et engendre le centre de $G$.
\subsubsection{}
Posons $a=x_{1}^{x_{2}}$, $b=x_{3}$, $c=x_{4}$, $d=x_{5}$, $e=x_{2}^{x_{3}x_{4}x_{5}}$, $f=x_{0}$; alors $a,b,c,d,e,f$ engendrent $G$ et satisfont à:
\[
\begin{picture}(150,140)

\thinlines             \put(52,48){$d$}
            \put(-50,72){$Q_{111}$}
             \put(20,72){$c$}
     \put(64,48){\line(1,0){25}}
     \put(94,48){$e$}
      
         \put(52,96){$b$} 
      \put(94,96){$a$}
           \put(64,98){\line(1,0){25}}
                  \put(104,96){\line(6,-5){25}}
        \put(104,51){\line(6,5){25}}
        \put(133,72){$f$}
  \put(27,72){\line(6,-5){25}}
              \put(27,72){\line(6,5){25}}
             \put(-50,25){$(f.a^{bcde})^{3}=1$}
              \end{picture}
\]
on a la relation $(adbecf)^{4}=z.(f^{ed}.a^{bc})^{3}$.
\subsubsection{}
Le groupe $G/Z(G)$, noté aussi $H^{*}_{3,6}$ admet la présentation
\[
 (a,b,\ldots,f/Q_{111},V=1)
 \]
  où $V=1$ est la relation hexagonale avec $V=(adbecf)^{4}$.\\ (\cite{[8]}, \cite{VD1}, \cite{VD3}, \cite{Z1}).
  \subsection{Le graphe $Q_{211}$}
  On désigne par $G$ le groupe $G^{+}(6,3)$ (noté $2.O_{6}^{-}(3):2$ dans l'ATLAS) et par $\tilde{G}$ l'extension centrale non scindée
  \[
  1\longrightarrow 3 \longrightarrow \tilde{G} \longrightarrow  G \rightarrow 1.
  \]
\subsubsection{}
Le groupe $\tilde{G}$ admet la présentation $(x_{i} (1\leqslant i \leqslant 6)/g_{6},(x_{3}^{x_{4}}x_{5})^{3}=1$) avec $g_{6}$:
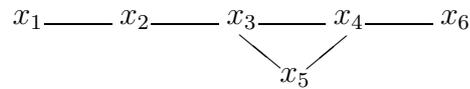
\[
\begin{picture}(150,100)
\thinlines
\put(-30,72){$x_{1}$}
\put(-18,72){\line(1,0){25}}
\put(10,72){$x_{2}$}
\put(22,72){\line(1,0){25}}
\put(50,72){$x_{3}$}
\put(62,72){\line(1,0){25}}
\put(90,72){$x_{4}$}
\put(102,72){\line(1,0){25}}
\put(130,72){$x_{6}$}
\put(55,68){\line(6,-5){15}}
\put(70,50){$x_{5}$}
\put(80,55){\line(6,5){15}}
\end{picture}
\]
Il existe un unique élément $x_{0}$ tel que:
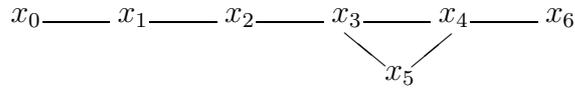
\[
\begin{picture}(150,100)
\thinlines
\put(-70,72){$x_{0}$}
\put(-58,72){\line(1,0){25}}
\put(-30,72){$x_{1}$}
\put(-18,72){\line(1,0){25}}
\put(10,72){$x_{2}$}
\put(22,72){\line(1,0){25}}
\put(50,72){$x_{3}$}
\put(62,72){\line(1,0){25}}
\put(90,72){$x_{4}$}
\put(102,72){\line(1,0){25}}
\put(130,72){$x_{6}$}
\put(55,68){\line(6,-5){15}}
\put(70,50){$x_{5}$}
\put(80,55){\line(6,5){15}}
\end{picture}
\]
On a: $x_{0}=x_{1}^{xx'x_{2}x_{1}}$, $x=x_{3}^{x_{4}x_{2}x_{5}x_{6}x_{3}x_{4}x_{3}}$, $x'=x_{3}^{xx_{2}x_{5}x_{3}}$.
Le centre de $\tilde{G}$ est d'ordre 6, $Z(\tilde{G})=<m,z>$. L'élément $z$, d'ordre 3, engendre le centre du sous-groupe $<x_{i} (0 \leqslant i \leqslant 5)>$ isomorphe à $H_{3,6}$ (voir l'Annexe 1). L'élément $m$ (produit de six involutions commutant deux à deux) est d'ordre 2 et engendre le centre du sous-groupe $<x_{3}^{x_{4}},x_{i} (0 \leqslant i \leqslant 6,i\ne 3,4)>$ isomorphe à $W(D_{6})$. On a 
\[
m=x_{2}x_{5}x_{6}xx'x_{0}=(x_{3}^{x_{4}}x_{6}x_{5}x_{2}x_{1}x_{0})^{5}.
\]
\subsubsection{}
Posons $a=x_{3}^{x_{4}}$, $b=x_{2}$, $c=x_{1}$, $d=x_{0}$, $e=x_{3}^{x_{2}x_{1}x_{0}}$, $f=x_{5}$, $a'=x_{6}$. Alors les éléments $a,b,\ldots,f,a'$ engendrent $\tilde{G}$ et satisfont à
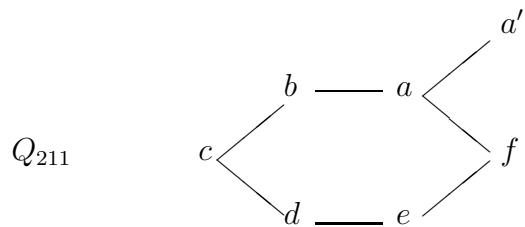
\[
\begin{picture}(150,140)

\thinlines             \put(52,48){$d$}
            \put(-50,72){$Q_{211}$}
             \put(20,72){$c$}
     \put(64,48){\line(1,0){25}}
     \put(94,48){$e$}
      
         \put(52,96){$b$} 
      \put(94,96){$a$}
           \put(64,98){\line(1,0){25}}
                  \put(104,96){\line(6,-5){25}}
                  \put(104,96){\line(6,5){25}}
                  \put(133,120){$a'$}
        \put(104,51){\line(6,5){25}}
        \put(133,72){$f$}
  \put(27,72){\line(6,-5){25}}
              \put(27,72){\line(6,5){25}}
              \end{picture}
\]
\subsubsection{}
Le groupe $G$ ($G \simeq G^{+}(6,3)$) admet la présentation:
\[ 
G=(a,b,\ldots,f,a'/Q_{211},V=1)
\] 
où $V=1$ est la relation hexagonale ($V=(adbecf)^{4}$); le centre de $G$ est engendré par l'involution $m$, $m$ est l'involution centrale de tout sous-groupe isomorphe à $W(D_{6})$ engendré par des conjugués de $a$ dans $G$. Posons $\mathcal{A}=a^{a'bfa}$; on a les égalités suivantes:
\[
m=a'bfd\mathcal{A}\mathcal{A}^{cbdc}=(aa'bfcd)^{5}=a'bfd\mathcal{A}\mathcal{A}^{edfe}=(aa'bfed)^{5}
\]
(\cite{[2]}, \cite{[8]}, \cite{VD1}, \cite{VD3}, \cite{Z1}).
\subsection{Le graphe $Y_{321}$}
On désigne par $G$ le groupe de Coxeter $W(E_{7})$ ($\simeq 2\times Sp_{6}(2)\simeq 2\times O_{7}(2)$)
\subsubsection{}
Le groupe $G$ admet la présentation $(a,b,\ldots,e,c',e'/Y_{321})$ avec:
\[
\begin{picture}(150,140)

\thinlines  
\put(-18,72){$c'$}    
\put(-10,72){\line(1,0){25}}
       \put(52,48){$d$}
            \put(-80,72){$Y_{321}$}
             \put(20,72){$c$}
     \put(64,48){\line(1,0){25}}
     \put(94,48){$e$}
      
         \put(52,96){$b$} 
      \put(94,96){$a$}
           \put(64,98){\line(1,0){25}}
                 \put(190,72){;}
                  \put(104,98){\line(1,0){25}}
                 \put(133,96){($a^{o}$)}
        \put(104,48){\line(1,0){25}}
        \put(133,48){$e'$}
  \put(27,72){\line(6,-5){25}}
              \put(27,72){\line(6,5){25}}
              \end{picture}
\]
il existe dans $G$ un unique élément $a^{o}$ qui complète le diagramme; $a^{o}$ s'écrit:\\
$a^{o}=\mathcal{C}^{abedcc'e'edcba}$ avec $\mathcal{C}=c^{c'bdc}$.
\subsubsection{}
Le centre de $G$ est engendré par une involution $z_{1}$ (produit de sept éléments d'ordre 2 commutant deux à deux); on a: 
\[
z_{1}=c'bd\mathcal{C}\mathcal{C}^{ee'de}e'a^{o}=c'bd\mathcal{C}\mathcal{C}^{aa^{o}ba}e'a^{o}=(cc'bdeae')^{9}=(cc'bdeaa^{o})^{9}
\]
\subsubsection{}(Avec les notations ci-dessus). Le groupe $\tilde{G}$ avec la présentation 
\[
(a,b,\ldots,e,a',c',e'/Y_{331},[z_{1},a']=1)
\]
\[
\begin{picture}(150,140)

\thinlines  
\put(-18,72){$c'$}    
\put(-10,72){\line(1,0){25}}
       \put(52,48){$d$}
            \put(-80,72){$Y_{331}$}
             \put(20,72){$c$}
     \put(64,48){\line(1,0){25}}
     \put(94,48){$e$}
      
         \put(52,96){$b$} 
      \put(94,96){$a$}
           \put(64,98){\line(1,0){25}}
                
                  \put(104,98){\line(1,0){25}}
                   \put(133,96){$a'$}
        \put(104,48){\line(1,0){25}}
        \put(133,48){$e'$}
  \put(27,72){\line(6,-5){25}}
              \put(27,72){\line(6,5){25}}
              \end{picture}
\]
est isomorphe à $A \rtimes W(E_{7})$, où $A$ est un 2-groupe abélien élémentaire d'ordre $2^{7}$ engendré par les conjugués de $a^{o}a'$; le centre de $\tilde{G}$ est engendré par $z_{1}$ et $z_{2}$ ($z_{2}$ désignant l'involution centrale du sous-groupe $<a,b,\ldots,e,a',c'>$), $z_{2}=(cc'bdaea')^{9}$.
Le produit $z_{1}z_{2}$ est dans $A$ et les relations $[z_{1},a']=1$ et $(a^{o}a')^{2}=1$ sont équivalentes.\\ (\cite{[1]}, \cite{[2]}, \cite{[3]}, \cite{[5]}, \cite{[7]}, \cite{VD2}, \cite{VD4}, \cite{Z2})
\subsection{Le graphe $Q_{221}$}
On désigne par $U$ le groupe simple $U_{6}(2)$ ($\simeq$  $^{2}A_{5}(2)$), son multiplicateur de Schur est $2^{2}\times 3$. Les éléments $t_{i}$ qui interviennent dans cette section sont des transvections unitaires et n'ont pas de lien avec les $t_{i}$ introduits section 3.
\subsubsection{}
Le groupe $\hat{K}$ isomorphe à $(2^{2}\times 3).2\times U$ admet la présentation
\[
(t_{i} (1 \leqslant i \leqslant 6/\underline{u},(t_{13}t_{6})^{2}=1)\qquad \mbox{avec}\qquad \underline{u}:
\]
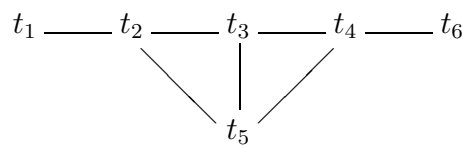
\[
\begin{picture}(150,100)
\thinlines
\put(-30,72){$t_{1}$}
\put(-18,72){\line(1,0){25}}
\put(10,72){$t_{2}$}
\put(22,72){\line(1,0){25}}
\put(50,72){$t_{3}$}
\put(62,72){\line(1,0){25}}
\put(90,72){$t_{4}$}
\put(102,72){\line(1,0){25}}
\put(130,72){$t_{6}$}
\put(55,68){\line(0,-1){25}}
\put(50,32){$t_{5}$}
\put(62,38){\line(1,1){28}}
\put(18,66){\line(1,-1){28}}
\put(-30,10){$(t_{3}^{t_{2}}t_{5})^{3}=(t_{3}^{t_{4}}t_{5})^{3}=(t_{3}^{t_{2}t_{4}}t_{5})^{3}=1,$}
\end{picture}
\]
où $t_{13}=t_{1}^{t't''t_{3}t_{4}}$ ($t'=t_{4}^{t_{3}t_{5}t_{2}t_{3}t_{5}}$, $t''=t_{2}^{t_{3}t_{5}t_{4}t_{3}t_{5}}$).
Les éléments $t_{i}$ correspondent aux transvections unitaires $t_{v_{i}}$, les $v_{i}$ ($1 \leqslant i \leqslant 6$) forment une base de l'espace vectoriel sur $\mathbb{F}_{4}$ associé et satisfont à $(v_{i},v_{j})=0$ si $(t_{i}t_{j})^{2}=1$, $(v_{i},v_{j})=1$ si $(t_{i}t_{j})^{3}=1$, sauf pour $\{i,j\}=\{3,5\}$, où on a alors $(v_{5},v_{3})=\omega$ avec $\mathbb{F}_{4}=\{0,1,\omega,\omega^{2}\}$.
L'élément $t_{13}$ correspond à la transvection $t_{v_{1}+v_{3}}$ et satisfait aux relations:
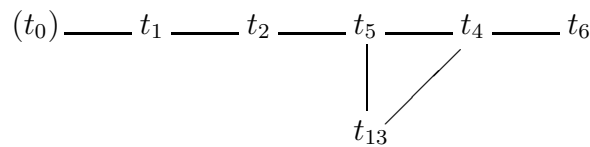
\[
\begin{picture}(150,100)
\thinlines
\put(-78,72){($t_{0}$)}
\put(-58,72){\line(1,0){25}}
\put(-30,72){$t_{1}$}
\put(-18,72){\line(1,0){25}}
\put(10,72){$t_{2}$}
\put(22,72){\line(1,0){25}}
\put(50,72){$t_{5}$}
\put(62,72){\line(1,0){25}}
\put(90,72){$t_{4}$}
\put(102,72){\line(1,0){25}}
\put(130,72){$t_{6}$}
\put(55,68){\line(0,-1){25}}
\put(50,32){$t_{13}$}
\put(62,38){\line(1,1){28}}
\put(30,10){$(t_{5}^{t_{4}}.t_{13})^{3}=1.$}
\end{picture}
\]
Le sous-groupe $T=<t_{13},t_{i} (1 \leqslant i \leqslant 6, i\ne 3)>$ est isomorphe à un quotient du groupe $\tilde{G}$ (voir Annexe 2), son centre contient un élément d'ordre 3, central dans $\hat{K}$.
\subsubsection{}
Il existe des conjugués de $t_{1}$ dans $\hat{K}$: $a,b,\ldots,f,a',c'$ tels que:
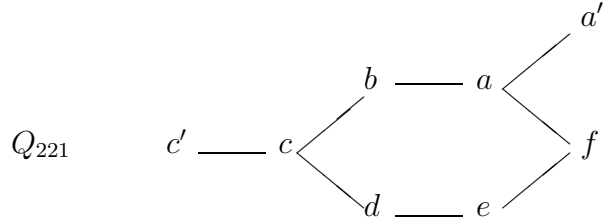
\[
\begin{picture}(150,140)

\thinlines 
\put(-22,72){$c'$}
\put(-10,72){\line(1,0){25}}
\put(52,48){$d$}
\put(-80,72){$Q_{221}$}
\put(20,72){$c$}
\put(64,48){\line(1,0){25}}
\put(94,48){$e$}
\put(52,96){$b$} 
\put(94,96){$a$}
\put(64,98){\line(1,0){25}}
\put(104,96){\line(6,-5){25}}
\put(104,96){\line(6,5){25}}
\put(133,120){$a'$}
\put(104,51){\line(6,5){25}}
\put(133,72){$f$}
\put(27,72){\line(6,-5){25}}
\put(27,72){\line(6,5){25}}
\end{picture}
\]
On a la correspondance:
\begin{eqnarray*}
a \longrightarrow t_{v_{4}+v_{5}}, b \longrightarrow t_{v_{2}}, c \longrightarrow t_{v_{1}}, d \longrightarrow t_{v_{2}+\omega (v_{1}+v_{3})+\omega v_{6}},\\ e \longrightarrow t_{\omega^{2}v_{1}+\omega(v_{1}+v_{6})+v_{5}}, f \longrightarrow t_{v_{1}+v_{3}}, a' \longrightarrow t_{v_{6}}, c' \longrightarrow t_{v_{2}+v_{6}}.
\end{eqnarray*}

\subsubsection{}
Le groupe $K$ (isomorphe à $2^{2}.2.U$) admet la présentation 
\[
(a,b,\ldots,f,a',c'/Q_{221}, V=1) 
\]
où $V=1$ désigne la relation hexagonale $V=(adbecf)^{4}$.

a) 
Il existe un unique élément $e^{o}$ satisfaisant aux relations:
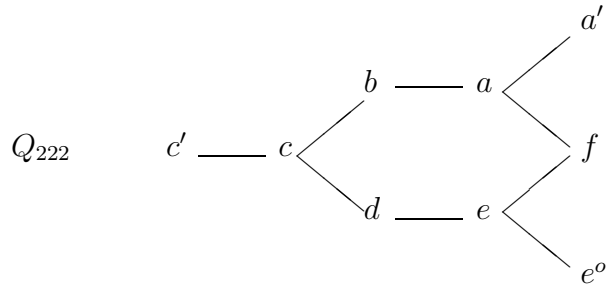
\[
\begin{picture}(150,140)

\thinlines 
\put(-22,72){$c'$}
\put(-10,72){\line(1,0){25}}
\put(52,48){$d$}
\put(-80,72){$Q_{222}$}
\put(20,72){$c$}
\put(64,48){\line(1,0){25}}
\put(94,48){$e$}
\put(52,96){$b$} 
\put(94,96){$a$}
\put(64,98){\line(1,0){25}}
\put(104,96){\line(6,-5){25}}
\put(104,96){\line(6,5){25}}
\put(133,120){$a'$}
\put(104,51){\line(6,5){25}}
\put(133,72){$f$}
\put(27,72){\line(6,-5){25}}
\put(27,72){\line(6,5){25}}
\put(104,51){\line(6,-5){25}}
\put(133,24){$e^{o}$}
\end{picture}
\]
cet élément s'écrit $ e^{o}=\mathcal{C}^{abedcc'a'abcde}=\mathcal{A}^{cbedaa'c'cbafe}$
(avec $\mathcal{A}=a^{a'bfa}$ et $\mathcal{C}=c^{c'bdc}$) et correspond à la transvection $t_{v_{1}+v_{3}+\omega v_{6}}$.

b) Les sous-groupes $R_{a}$, $R_{c}$ et $R_{e}$ définis à partir des graphes:
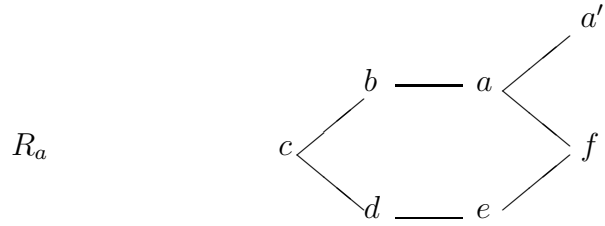
\[
\begin{picture}(150,140)

\thinlines 
\put(52,48){$d$}
\put(-80,72){$R_{a}$}
\put(20,72){$c$}
\put(64,48){\line(1,0){25}}
\put(94,48){$e$}
\put(52,96){$b$} 
\put(94,96){$a$}
\put(64,98){\line(1,0){25}}
\put(104,96){\line(6,-5){25}}
\put(104,96){\line(6,5){25}}
\put(133,120){$a'$}
\put(104,51){\line(6,5){25}}
\put(133,72){$f$}
\put(27,72){\line(6,-5){25}}
\put(27,72){\line(6,5){25}}
\end{picture}
\]
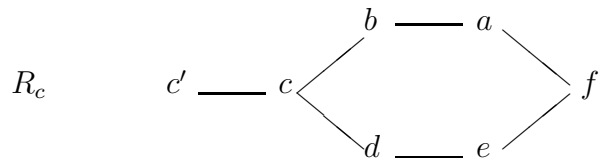
\[
\begin{picture}(150,90)

\thinlines 
\put(-22,72){$c'$}
\put(-10,72){\line(1,0){25}}
\put(52,48){$d$}
\put(-80,72){$R_{c}$}
\put(20,72){$c$}
\put(64,48){\line(1,0){25}}
\put(94,48){$e$}
\put(52,96){$b$} 
\put(94,96){$a$}
\put(64,98){\line(1,0){25}}
\put(104,96){\line(6,-5){25}}
\put(104,51){\line(6,5){25}}
\put(133,72){$f$}
\put(27,72){\line(6,-5){25}}
\put(27,72){\line(6,5){25}}
\end{picture}
\]
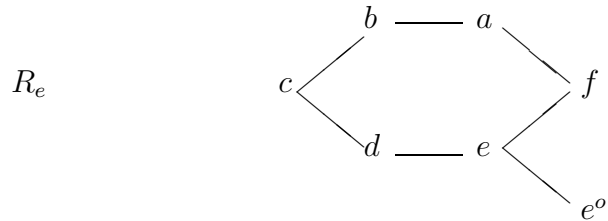
\[
\begin{picture}(150,90)

\thinlines 
\put(52,48){$d$}
\put(-80,72){$R_{e}$}
\put(20,72){$c$}
\put(64,48){\line(1,0){25}}
\put(94,48){$e$}
\put(52,96){$b$} 
\put(94,96){$a$}
\put(64,98){\line(1,0){25}}
\put(104,96){\line(6,-5){25}}
\put(104,51){\line(6,5){25}}
\put(133,72){$f$}
\put(27,72){\line(6,-5){25}}
\put(27,72){\line(6,5){25}}
\put(104,51){\line(6,-5){25}}
\put(133,24){$e^{o}$}
\end{picture}
\]
sont isomorphes au groupe orthogonal $2.O_{6}^{-}(3):2$ (voir Annexe 2); leurs images dans $U$ sont les représentants des trois classes de conjugaison de sous-groupes isomorphes à $O_{6}^{-}(3):2$ . Leurs involutions centrales $\mu_{a}$, $\mu_{c}$, $\mu_{e}$ appartiennent au multiplicateur de Schur de $U$, on a $\mu_{a}\mu_{c}\mu_{e}=1$, $\mu_{a}=(aa'fbcd)^{5}$, $\mu_{c}=(cc'dbaf)^{5}$ et $\mu_{e}=(ee^{o}dfab)^{5}$.
\newpage
c) Les sous-groupes $W_{a}$, $W_{c}$ et $W_{e}$ définis à partir des graphes:
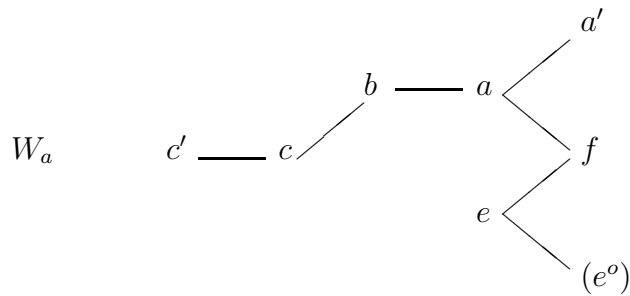
\[
\begin{picture}(150,140)

\thinlines 
\put(-22,72){$c'$}
\put(-10,72){\line(1,0){25}}
\put(-80,72){$W_{a}$}
\put(20,72){$c$}
\put(94,48){$e$}
\put(52,96){$b$} 
\put(94,96){$a$}
\put(64,98){\line(1,0){25}}
\put(104,96){\line(6,-5){25}}
\put(104,96){\line(6,5){25}}
\put(133,120){$a'$}
\put(104,51){\line(6,5){25}}
\put(133,72){$f$}
\put(27,72){\line(6,5){25}}
\put(104,51){\line(6,-5){25}}
\put(133,24){$(e^{o})$}
\end{picture}
\]
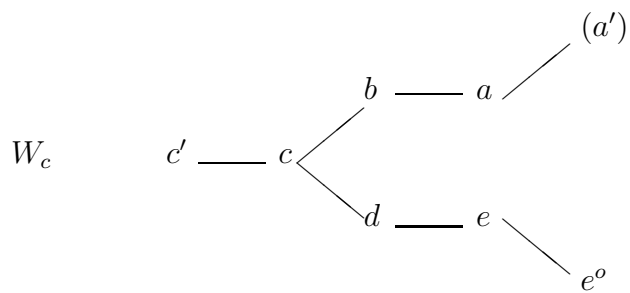
\[
\begin{picture}(150,140)

\thinlines 
\put(-22,72){$c'$}
\put(-10,72){\line(1,0){25}}
\put(52,48){$d$}
\put(-80,72){$W_{c}$}
\put(20,72){$c$}
\put(64,48){\line(1,0){25}}
\put(94,48){$e$}
\put(52,96){$b$} 
\put(94,96){$a$}
\put(64,98){\line(1,0){25}}
\put(104,96){\line(6,5){25}}
\put(133,120){$(a')$}
\put(27,72){\line(6,-5){25}}
\put(27,72){\line(6,5){25}}
\put(104,51){\line(6,-5){25}}
\put(133,24){$e^{o}$}
\end{picture}
\]
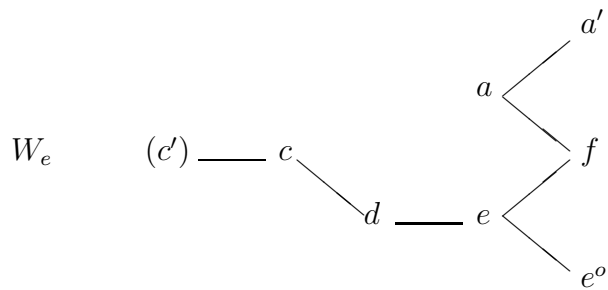
\[
\begin{picture}(150,140)

\thinlines 
\put(-30 ,72){$(c')$}
\put(-10,72){\line(1,0){25}}
\put(52,48){$d$}
\put(-80,72){$W_{e}$}
\put(20,72){$c$}
\put(64,48){\line(1,0){25}}
\put(94,48){$e$}
\put(94,96){$a$}
\put(104,96){\line(6,-5){25}}
\put(104,96){\line(6,5){25}}
\put(133,120){$a'$}
\put(104,51){\line(6,5){25}}
\put(133,72){$f$}
\put(27,72){\line(6,-5){25}}
\put(104,51){\line(6,-5){25}}
\put(133,24){$e^{o}$}
\end{picture}
\]
contiennent respectivement $e^{o}$,$a'$ et $c'$; ils sont isomorphes à $W(E_{7})$, leurs involutions centrales $z_{a}$, $z_{c}$ et $z_{e}$ sont dans $Z(K)$. 

On a $\mu_{a}z_{a}=\mu_{c}z_{c}=\mu_{e}z_{e}=q$ où $q$ correspond au produit des cinq transvections d'un plan isotrope. On a $Z(K)=<\mu_{a},\mu_{c}>.<\mu_{a}z_{a}>$.\\ (\cite{[2]}, \cite{[4]}, \cite{VD1}, \cite{VD2},\cite{VD3}, \cite{VD5}, \cite{Z3}).

\bibliographystyle{plain}
\bibliography{MM1}

\end{document}